# ASYMPTOTIC BEHAVIOR OF SYSTEMS INVOLVING DELAYS: PRELIMINARY RESULTS

M. De la Sen

Instituto de Investigación y Desarrollo de Procesos IIDP

Departamento de Electrticidad y Electrónica. **Facultad de Ciencia y Tecnología**

Leioa (Bizkaia). Aptdo. 644 de Bilbao. 48080- Bilbao. SPAIN

**Abstract**. This paper investigates the relations between the particular eigensolutions of a limiting functional differential equation of any order, which is the nominal (unperturbed) linear autonomous differential equations, and the associate ones of the corresponding perturbed functional differential equation. Both differential equations involve point and distributed delayed dynamics. The proofs are based on a Perron type theorem for functional equations so that the comparison is governed by the real part of a dominant zero of the characteristic equation of the nominal differential equation. The obtained results are also applied to investigate the global stability of the perturbed equation based on that of its corresponding limiting equation.

**Keywords**. Time- delay systems, Perron type theorem, Point delays, Distributed delays, Comparison of solutions, Asymptotic behavior.

**1. Introduction**

Time-delay dynamic systems are an interesting field of research in dynamic systems and functional differential equations because of intrinsic theoretical interest because the formalism lies in that of functional differential equations, then infinite dimensional and because of the wide range of applicability in modelling of physical systems, like transportation systems, queuing systems, tele-operated systems, war/ peace models, Biological systems, finite impulse response filtering, etc., [8], [11], [18], [20]. Important particular interest has been devoted to stability, stabilization and model-matching of control systems where the object to be controlled possess delayed dynamics and the controller is synthesized either with delayed dynamics or it is delay-free (see, for instance, [3-5], [8-10], [12-15], [17], [20]). The properties are formulated as either being independent of or dependent on the sizes of the delays. An intrinsic problem which generated analysis complexity is the presence of infinitely may characteristic zeros because of the functional nature of the dynamics. This fact generates difficulties in the closed-loop pole–placement problem compared to the delay-free case, [17], as well as in the stabilization problem, [3-4], [9-13], [15], [20], [23-28], including the case of singular time-delay systems where the solution is sometimes non-unique, and impulsive, because of the dynamics associated to a nilpotent matrix, [23]. The properties of the associated evolution operators have been investigated in [4], [11] and [13]. This paper is devoted to obtain results relying on a comparison and an asymptotic comparison of the eigensolutions between a nominal (unperturbed) functional differential equation involving wide classes of delays and a perturbed version (describing the current dynamics) with some smallness in the limit assumptions on the perturbed functional differential equation. The nominal equation is the limiting equation of the perturbed one since the parameters of the last one converge asymptotically to those of its limiting counterpart. The problem is of interest in practice is since very often the perturbations related to a nominal model in dynamic systems occur during the transients while they are asymptotically vanishing in the steady-state or, in the more general worst case, they grow at a smaller rate than the solution of the nominal differential equation. In this context, the nominal differential equation may be viewed as the limiting equation of the perturbed one. The comparison between the solutions of the limiting differential equation and those of the perturbed one based on Perron-type results have been studied classically for ordinary differential equations and more recently for the case of functional equations [1-2], [19]. Particular functional equations of interest are those involving both point and distributed delays potentially including the last ones Volterra-type terms, [3-6], [11].

**Notation**: $R_{0+} := R_+ \cup \{0\}$, $R_+ := \{z \in R: z > 0\}$, $R_{0-} := R_- \cup \{0\}$, $R_- := \{z \in R: z < 0\}$,
$C_{0+} = \{z \in C: Re\, z \geq 0\}$, $C_+ := \{z \in C: Re\, z > 0\}$, $C_{0-} := \{z \in C: Re\, z \leq 0\}$, $C_- := \{z \in C: Re\, z < 0\}$
$Z_{0+} := Z_+ \cup \{0\}$, $Z_+ := \{z \in Z: z > 0\}$

where **R**, **C** and **Z** are the sets of real, complex and integer numbers, respectively. The complex imaginary unity is $i = \sqrt{-1}$. A finite subset of j consecutive positive integers starting with 1 is denoted by $\bar{j} := \{1, 2, ..., j\}$. The set $R_{-h} = [-h, 0) \cup R_{0+}$ will be used to define the solution of functional differential equations on $R_{0+}$ including its initial condition on $[-h, 0]$.
$C^{(i)}(R_{0+}, R_{0+}^m)$ is the set of m-vector real functions of class $C^{(i)}$ and definition domain $R_{0+}$ and $PC^{(i)}(R_{0+}, R_{0+}^m)$ is the set of m-vector real functions in $C^{(i-1)}(R_{0+}, R_{0+}^m)$ whose i-th derivative is piecewise continuous. Similar sets of functions are defined when the ranges are complex as $C^{(i)}(R_{0+}, C_{0+}^m)$ and $PC^{(i)}(R_{0+}, C_{0+}^m)$.

For the delayed system, $T: [0, \infty) \to L(X)$ is the inverse Laplace transform of the resolvent mapping $\hat{T}(s)$, which is holomorphic where it exists, with X being the real Banach space of n-vector real functions endowed with the supremum norm on their definition domain defined for any such a complex or real vector function $\phi$ of definition domain $D$ by $|\phi|_\alpha = \sup_{\tau \in D} (\|\phi(\tau)\|_\alpha)$ where $\|.\|_\alpha$ denotes any of the standard vector norms for $\alpha = 2$; i.e. Euclidean (also called Froebenius or $\ell_2$)- norm, $\alpha = 1$ for the $\ell_1 - norm$, $\alpha = \infty$ for the $\ell_\infty - norm$, etc. Similar notations are used for the corresponding matrix induced norms. The un-subscripted symbol $|.|$ is used for absolute values of real, complex and integer numbers, as usual. It is said that the delays associated with Volterra-type dynamics are infinitely distributed because of the contribution of the delayed dynamics is made under an integral over $[0, \infty)$ as $t \to \infty$, i.e. $x(t - \tau - h_i')$ acts on the dynamics of $x(t)$ from $\tau = 0$ to $\tau = t$ for finite t and as $t \to \infty$.

Dom (H) is the definition domain of the operator H and $sp(A)$ is the spectrum (i.e. the set of distinct eigenvalues) of the square matrix A. The matrix measure of the norm- dependent complex-valued matrix A is defined by $\kappa_\alpha(A) := \lim\limits_{\delta \to 0^+} \dfrac{\|I_n + \delta A\|_\alpha - \delta}{\delta} \geq \operatorname{Re}\lambda_i(A), \ \forall \lambda_i \in sp(A)$.

## 2. Problem statement and basic first results

Consider the following linear nominal functional differential systems with point and, in general, both Volterra- type and finite distributed delays:

$$\dot{x}(t) = Lx_t \equiv \sum_{i=0}^{m} A_i x(t - h_i) + \sum_{i=0}^{m'} \int_0^t d\alpha_i(\tau) A_{\alpha_i} x(t - \tau - h_i') + \sum_{i=m'+1}^{m'+m''} \int_{t-h_i'}^{t} d\alpha_i(t-\tau) A_{\alpha_i} x(\tau)$$

(2.1)

and

$$\dot{x}(t) = Lx_t + f(t, x_t) \equiv \sum_{i=0}^{m} A_i x(t - h_i) + \sum_{i=0}^{m'} \int_0^t d\alpha_i(\tau) A_{\alpha_i} x(t - \tau - h_i') + \sum_{i=m'+1}^{m'+m''} \int_{t-h_i'}^{t} d\alpha_i(t-\tau) A_{\alpha_i} x(\tau) + f(t, x_t)$$

(2.2)

$$f(t, x_t) = \sum_{i=0}^{m} \tilde{A}_i(t) x(t - h_i) + \sum_{i=0}^{m'} \int_0^t d\alpha_i(\tau) \tilde{A}_{\alpha_i}(t) x(t - \tau - h_i') + \sum_{i=m'+1}^{m'+m''} \int_{t-h_i'}^{t} d\alpha_i(t-\tau) \tilde{A}_{\alpha_i}(\tau) x(\tau) + f_0(t, x_t)$$

(2.3)

Eq. (2.1) is the limiting equation of the perturbed equation (2.2), subject to (2.3), for $f(t, x_t) \to 0$ as $t \to \infty$ under the following technical hypothesis:

**H.1**: The initial conditions of both differential equations (2.1) and (2.2) are real n-vector functions $\phi \in C_e(-h)$ where $C_e(-h) := \{\phi = \phi_1 + \phi_2 : \phi_1 \in C(-h), \phi_2 \in B^0(-h)\}$, $\phi(0) = x_0$, with $C(-h) := \{C^0([-h, 0], X)\}$; i.e., the set of continuous mappings from $[-h, 0]$ into the Banach space X with norm $\bar{\phi}_\alpha := |\phi|_\alpha = \operatorname{Sup}\{\|\phi(t)\|_\alpha : -h \leq t \leq 0\}$ ; $\|\cdot\|$ denoting the Euclidean norm of vectors in $C^n$ and matrices in $C^{n \times n}$, and $B^0(-h) := \{\phi : [-h, 0] \to X\}$ is the set of real bounded vector functions on X endowed with the supremum norm having support of zero measure. Roughly speaking, $\phi \in B^0(-h)$ if and only if it is almost everywhere zero except at isolated discontinuity points within $[-h, 0]$ where it is bounded. Thus, $\phi \in C_e(-h)$ if and only if it is almost everywhere continuous in $[-h, 0]$ except possibly on a set of zero measure of bounded discontinuities. $C_e(-h)$ is also endowed with the supremum norm since $\phi = \phi_1 + \phi_2$, some $\phi_1 \in C(-h)$, $\phi_2 \in B^0(-h)$ for each $\phi \in C_e(-h)$. In the following, the supremum norms on $L(X)$ are also denoted with $|\cdot|$.

Close spaces of functions are $C(\boldsymbol{R}_{-h}) := C([-h, \infty), C^n)$ which is the Banach space of continuous functions from $[-h, \infty)$ into $C^n$ endowed with the norm $|\phi|_\alpha = \sup\limits_{-h \leq \tau < \infty}(\|\phi(\tau)\|_\alpha)$; $\forall \phi \in C_e(-h) := C([-h, 0), C^n)$ being an initial condition, for some given vector norm $\|\cdot\|_\alpha$, . Note that for $t \in \boldsymbol{R}_{0+}$, the solution which satisfies (2.2) subject to (2.3) is in $C_e(\boldsymbol{R}_{0+}) := C(\boldsymbol{R}_{0+}, C^n)$, the Banach space of continuous functions from $\boldsymbol{R}_{0+}$ into $C^n$ which satisfies (2.2)-(2.3), $\forall \phi \in C_e(-h)$ , endowed with $|\phi|_\alpha = \sup\limits_{0 \leq \tau < \infty}(\|\phi(\tau)\|_\alpha)$ .

Thus, $L: C_e(\mathbf{R}_{-h}) \to \mathbf{C}^n$ is a bounded linear functional defined by the right hand-side of (2.1).

**H.2**: All the operators $A_k (0 \leq k \leq m)$, $A_{\alpha_k} (0 \leq k \leq m' + m'')$ are in $L(X) := L(X, X)$, the set of linear operators on $X$, of dual $X^*$, and $h_k$ and $h'_\ell$ ($k = 1, 2, \ldots, m$; $\ell = 0, 1, \ldots, m' + m''$) are nonnegative constants with $h_0 = h'_0 = 0$ and $h := \mathrm{Max}\left( \underset{1 \leq i \leq n}{\mathrm{Max}}(h_i), \underset{1 \leq i \leq m' + m''}{\mathrm{Max}}(h'_i) \right)$.

**H3**: The linear operators $A_{\alpha_i} \in L(X)$, with abbreviated notation $A_{\alpha_0} = A_\alpha$, are closed and densely defined linear operators with respective domain and range $D(A_{\alpha_i})$ and $R(A_{\alpha_i}) \subset X$ ($i = 0, 1, \ldots, m' + m''$). The functions $\alpha_i \in C^0([0,\infty), \mathbf{C}) \cap BV_{loc}(\mathbf{C}_+)$ ($i = 0, 1, \ldots, m'$) and $\alpha_i \in C^0([-h, 0), \mathbf{C})$ (i=0, 1,..., m' + m") being every- where differentiable with possibly bounded discontinuities on subsets of zero measure of their definition domains with $\int_0^\infty e^{vt} |d\alpha_i(t)| < \infty$ some nonnegative real constant $v$ ($i = 0, 1, \ldots, m'$). If $\alpha_i(.)$ is a matrix function $\alpha_i : [0, \infty) \times X^* \to L(X, X^*)$ then it is in $C^0([0,\infty), \mathbf{C}^{n \times n}) \cap BV_{loc}(\mathbf{C}_+^{n \times n})$ with $\int_0^\infty e^{vt} |d\alpha_i(t)| < \infty$ and its entries being everywhere time-differentiable with possibly bounded discontinuities on a subset of zero measure of their definition domains.

**H.4** It is assumed that $f_0 : \mathbf{R}_{0+} \times C(\mathbf{R}_{-h}) \to \mathbf{C}^n$ and $h_0 = h'_0 = 0$,

$$x_t^{C(\mathbf{R}_{-h})} := \begin{cases} x : [-h, \tau) \to X, & \tau \leq t \\ 0, & \tau > t \end{cases}$$

, satisfying $x(t) = \phi(t)$, $\forall t \in [-h, 0]$, is a string of the solution of (2.2)-(2.3). Another strings of the solution trajectory of interest in this manuscript are $x_t^{C(\mathbf{R}_{0+})}$ which pointwise defined by x(t) within the interval $[0, t]$ and zero, otherwise, and subject to the constraint $x_t = \phi_t$ within $[\max(-h, t-h), t]$ for any real       and being zero outside this interval. Finally, $x_t$ denotes the solution string within $[t-h, t]$ pointwise defined by the solution x(t) to (2.2)-(2.3) for each $t \in R_{-h}$ being zero outside $[t-h, t]$ and subject to the constraint $x_t = \phi_t$ within $[\max(-h, t-h), t]$ for any real $t \leq h$ and being zero outside this interval. □

$A_i$ and $A_{\alpha_k}$ and $\tilde{A}_i : [0, \infty) \to \mathbf{C}^{n \times n}$ and $\tilde{A}_i : [0, \infty) \to \mathbf{C}^{n \times n}$ ($i = 0, 1, \ldots, m$; $k = 0\ 1, \ldots, m' + m''$) belong to the spaces of constant real matrices and real matrix functions, respectively. The last ones are also unbounded operators on a Banach space of n-vector real functions $x \in X$ endowed with the supremum norm where the vectors of point and distributed constant delays are:

$\hat{h} := (0, h_1, h_2, \ldots, h_m)$ and $\hat{h}' = \left( \hat{h}_1^{'T} \vdots \hat{h}_2^{'T} \right)^T := \left( 0, h'_1, h'_2, \ldots, h'_{m'} \vdots h'_{m'+1}, h'_{m'+2}, \ldots, h'_{m'+m''} \right)^T$, respectively, with $h_i \geq 0$ and $h'_k \geq 0$ ($i = 1, 2, \ldots, m' + m''$) being, respectively, point and distributed delays, with $h_0 = h'_0 = 0$, $A_0 \equiv A$, $A_{\alpha_0} \equiv A_\alpha$ and $\alpha_0(.) \equiv \alpha(.)$. The first $m'$ distributed delays are associated with Volterra- type dynamics. In other words, the infinitely distributed delays give contributions $\int_0^t d\alpha_i(\tau) A_{\alpha_i} x(t - \tau - h'_i)$ with finite real constants $h'_i$ with ($i = 1, 2, \ldots m'$) to $\dot{x}(t)$ which are point delays under the integral symbol. The functions $\alpha_i : [0, \infty) \to \mathbf{C}$ and $\alpha_k : [0, h'_k] \to \mathbf{C}$ are continuously differentiable real functions within their definition domains except possibly on sets of zero measure where the time-derivatives have bounded discontinuities. All or some of the $\alpha_i(.)$ may be alternatively matrix functions $\alpha_i : [0, t] \to C^{n \times n}$ ($i = 0, 1, ..., m'$) for $t \in \mathbf{R}_+$ and $\alpha_i : [0, h'_k] \to C^{n \times n}$ ($i = m'+1, m'+2..., m'+m''$) with $\alpha_i(0) = 0$; $i = 0, 1, ..., m'+m''$. On the other hand, the perturbation vector function $f(t, x_t)$ in (2.2), defined in (2.3), with respect to the limiting equation (2.1) is defined by the function $f : \mathbf{R}_{0+} \times C(-\mathbf{R}_h) \to C^{n \times n}$ which describes a perturbed dynamics associated with the delays plus a perturbation function $f_0 : \mathbf{R}_{0+} \times C(-\mathbf{R}_h) \to C^{n \times n}$ which is not included in the remaining terms of the function f in (2.3). Note that both the delayed differential systems (2.2)-(2.3) and its limiting version (2.1) are very general since it includes point-delayed dynamics, like, for instance in typical war/ peace models or the so-called Minorski's problem appearing when controlling the lateral dynamics of a ship, [11]. It also includes real constants $h'_i$ ( i=0, 1, ... , m '), with $h'_0 = 0$, associated with infinitely distributed delayed contributions to the dynamics through integrals, related to the ; i = 0 , 1 , ... , m '. Such delays are relevant, for instance, in viscoelastic fluids, electrodynamics and population growth [3], [8-9]. In particular, an integro-differential Volterra's type term is also included through $h'_0 = 0$. Apart from those delays, the action of finite distributed delays characterized by real constants $h'_i$ (i = 0, 1, ...., m ' + m '' ) is also included in the limiting equation (2.1) and in the equation (2.3) . That kind of delays is well-known, for instance, in econometric models related to production rate, [9].

The integrability of the $\alpha_i(.)$-functions (or matrix functions) on [ t - $h'_i$ , t ] ; m ´ + 1 follows since their definition domain is bounded. The technical hypothesis H1-H4 guarantee the existence and uniqueness of the solution in $C(\mathbf{R}_{0+}) := C(\mathbf{R}_{0+}, C^n)$ of the functional differential systems (2.1) and (2.2)-(2.3) for each given initial condition $\phi \in C_e(-h)$. Take Laplace transforms in (2.1) by using the convolution theorem and the relations $d\alpha(\tau) = \dot{\alpha}(\tau)d\tau$. It follows that $d\hat{\alpha}_i(s) = s\hat{\alpha}_i(s) - \alpha_i(0)$, where $\hat{f}(s) := Lap\, f(t)$ denotes the Laplace transform of f(t). Thus, the unique solutions of both the limiting equation (2.1) and that of (2.2) –(2.3) in $C(\mathbf{R}_{0+})$, subject to (2.3); $\forall t \in \mathbf{R}_{0+}$, for the same given initial conditions $\phi \in C_e(-h)$ are, respectively, defined by:

$$y(t) = T(t, 0)x(0^+) + \int_{-h}^{0} T(t, \tau)\phi(\tau)U(\tau)d\tau \qquad (2.4)$$

$$x(t) = T(t, 0)x(0^+) + \int_{-h}^{0} T(t, \tau)\phi(\tau)U(\tau)d\tau + \int_{0}^{t} T(t, \tau)f(\tau, x_\tau)d\tau \qquad (2.5)$$

where U (t) is the unit step (Heaviside) function and T(t, τ) is the evolution operator , [4],[7],[11] , of the linear equation (2.1) whose Laplace transform, everywhere it exist, is given by the resolvent:

$$\hat{T}(s) := Lap\, T(t,0) = \left[ s\left( I_n - \sum_{i=0}^{m'} \hat{\alpha}_i(s) A_{\alpha_i} e^{-h'_i s} - \sum_{i=m'+1}^{m'+m''} \hat{\alpha}_i(s) A_{\alpha_i}\left(1 - e^{-h'_i s}\right) \right) - \sum_{i=0}^{m} A_i e^{-h_i s} + \right.$$

$$\left. \alpha(0)A_\alpha + \sum_{i=0}^{m'} \alpha_i(0) A_{\alpha_i} e^{-h'_i s} + \sum_{i=m'+1}^{m'+m''} \alpha_i(0) A_{\alpha_i}\left(1 - e^{-h'_i s}\right) \right]^{-1} \qquad (2.6)$$

As usual, it is said though the manuscript that (2.1) is the limiting equation of (2.2)-(2.3) irrespective of the fact that $f(t, x_t)$ converges or not to zero as $t \to \infty$. The evolution operator is a convolution operator so that $T(t,\tau) = T(t-\tau, 0) = T(t-\tau)$ if the Volterra – type dynamics is zero or if the associate differentials in the Riemann- Stieljes integrals $d\alpha_i(t) = \chi_i dt$ with the $\chi_i$ being real constants. In this case, the limiting linear functional differential equation is, furthermore, time –invariant. Note that the limiting equation (2.1) is guaranteed to be globally exponentially uniformly stable if and only if $\hat{T}(s)$ exists within some region including properly the right-complex plane. In other words, if it is compact for Re s >- $\alpha_0$, for some r constant $\alpha_0 \in \mathbf{R}$ located to the right of all the real parts of all the zeros of $det\,\hat{T}^{-1}(s)$ (also often called the characteristic zeros of the limiting equation (2.2) or, simply, its eigenvalues), since then all the entries of its Laplace transform T(t) decay with exponential rate on $\mathbf{R}_{0+}$ for $\phi \in C_e(-h)$ and then $|x(t)|$ decays with exponential rate on $\mathbf{R}_+$. The main results addressed in [5] and [9-12] is the investigation of the global uniform exponential stability of (1). The stability of the limiting system (2.1) is investigated in [3-4] provided that any auxiliary system formed with any of the additive parts of the dynamics of (2.1) has such a property and provided that an impulsive- solution dependent input exists. The compactness of the relevant input- output and input-state operators under forcing external inputs and impulsive forcing terms is also investigated in [4]. The basic mathematical tool used is that in those papers is that the unique solution of the homogeneous (2.1) for each function of initial conditions $\phi \in C_e(-h)$ may be equivalently written in infinitely many cases by first rewriting (2.1) by considering different 'auxiliary' reference homogeneous systems plus additional terms considered as forcing actions. The objective of this paper is to compare the solutions (2.1) and (2.2), subject to (2.3), of the limiting and current functional differential equations (2.1) and (2.2) by using a Perron – type result using a similar technique as that used in [2]. The subsequent theorem is a generalization of a classical Perron type theorem for ordinary differential equations to Eq. (2.2), subject to Eq. (2.3) compared to Eq. (2.1) (see [1, Chapter IV, Theorem 5] and of Theorem 1.1 in [2] for functional differential equations which includes several kinds of delays such as point and distributed delays and Volterra-type dynamics with infinite delays. The results extends the perturbation term to include constant upper-bounding terms in the perturbation functional (2.3) and characteristic zeros of the limiting equation (2.2) (i.e. zeros $det\,\hat{T}^{-1}(s)$) of multiplicity greater than unity (being degenerated or non-degenerated) in the limiting dynamics defined by Eq. (2.1).

**Theorem 2.1**. Let x be a solution of (2.2), subject to (2.3), on $\mathbf{R}_{0+}$ subject to initial conditions $\phi \in C_e(-h)$ such that

$$|f(t, x_t)|_\alpha \leq \gamma_\alpha(t)|x_t|_\alpha + K_{0\alpha}, \quad t \in \mathbf{R}_{0+} \tag{2.7}$$

for some norm-dependent $K_{0\alpha} \in \mathbf{R}_{0+}$ where $\gamma_\alpha \in C^{(0)}(\mathbf{R}_{0+}, \mathbf{R}_{0+})$ is also norm-dependent and satisfies $\int_t^{t+1} \frac{t^{\vartheta_k - 1}}{\vartheta_k!} e^{\beta \sigma_k t} \gamma_\alpha(s) ds \to 0$ as $t \to \infty$ where $\beta = 0$ if there is no Volterra term in (2.1) and $\beta = 1$, otherwise and $\sigma_k$ are the real parts of the zeros of $det\,\hat{T}^{-1}(s)$ of the limiting equation (2.2) with respective multiplicities $\vartheta_k$. Then, the following properties hold:

(i) $\lim_{t \to \infty} \int_t^{t+1} g_\alpha(s) \left( \left| \frac{f(s, x_s)}{x_s} \right|_\alpha - K_{0\alpha} \right) ds = 0$ (2.8)

where $g_\alpha \in PC^{(0)}(\mathbf{R}_{0+}, \mathbf{R}_{0+})$ is an indicator function defined by

$$g_\alpha(t) = \begin{cases} 0 & \text{if } \gamma_{0\alpha}(t) = \left| \frac{f(t, x_t)}{x_t} \right|_\alpha - K_{0\alpha} \leq 0 \\ 1 & \text{otherwise} \end{cases}$$

(ii) The real numbers $\mu_k = \mu_k(x) = \lim_{t \to \infty} \frac{\log|x_t|}{t^k}$ exist, are norm-independent and finite, $\forall k \geq k_1$ and some integer $k_1 \geq 1$ with $\mu_k = 0$, $\forall k > k_1$ or $\lim_{t \to \infty} e^{bt} x(t) = 0$, $\forall b \in \mathbf{R}$. If $\mu_{k_1} \geq 0$ or if (2.7) holds with $K_{0\alpha} = 0$ then either $\mu_{k_1}$ it is the real part of a zero of $\det \hat{T}^{-1}(s)$, for which the resolvent $\hat{T}(s)$ trivially exists and it is bounded, or $\lim_{t \to \infty} e^{bt} x(t) = 0$, $\forall b \in \mathbf{R}$.

(iii) Assume that all the zeros of $\det \hat{T}^{-1}(s)$ have real negative parts and (2.7) holds only for some constants $K_{0\alpha} \in \mathbf{R}_+$. Then, either the limits $\mu_k = \mu_k(x) = \lim_{t \to \infty} \frac{\log|x_t|}{t^k} = 0$ exist, $\forall k \geq k_1$, some integer $k_1 \geq 1$, and furthermore $\mu_1$ is not trivially the real part of a characteristic zero of (2.2), or $\lim_{t \to \infty} e^{bt} x(t) = 0$, $\forall b \in \mathbf{R}$.

*Proof*: (i) From (2.7), $\gamma_\alpha(t) \geq \max\left(0, \left| \frac{f(t, x_t)}{x_t} \right|_\alpha - K_{0\alpha}\right) = g_\alpha(t) \left( \left| \frac{f(t, x_t)}{x_t} \right|_\alpha - K_{0\alpha} \right) \geq 0$ so that:

$$0 = \lim_{t \to \infty} \int_t^{t+1} \gamma_\alpha(s) ds \geq \limsup_{t \to \infty} \int_t^{t+1} g_\alpha(s) \left( \left| \frac{f(s, x_s)}{x_s} \right|_\alpha - K_{0\alpha} \right) ds \geq 0$$

$$\Rightarrow \exists \lim_{t \to \infty} \int_t^{t+1} g_\alpha(s) \left( \left| \frac{f(s, x_s)}{x_s} \right|_\alpha - K_{0\alpha} \right) ds = 0$$

and Property (i) has been proved.

(ii) From (2.7), $|f(t, x_t)|_\alpha = \gamma_\alpha(t) |x_t|_\alpha + K_{0\alpha} - \omega_\alpha(t)$, $t \in \mathbf{R}_{0+}$, some $\omega_\alpha \in PC^{(0)}(\mathbf{R}_{0+}, \mathbf{R}_{0+})$. Then, one gets from (2.5)-(2.6):

$$\|x(t)\|_\alpha \leq K_1(\alpha) \frac{t^{\nu-1}}{\nu!} \left[ e^{\mu t} \|x(0^+)\|_\alpha + \left| \frac{e^{\mu t}(1 - e^{\mu h})}{\mu} \right| \left( (|\phi|_\alpha + K_{0\alpha}) + \beta \max(1, e^{\mu t}) \sup_{0 \leq \tau \leq t} (|x_\tau|_\alpha) \left( \sum_{l=0}^{j_t} \int_\ell^{\ell+1} \gamma_\alpha(s) ds + \int_{j_t}^t \gamma_\alpha(s) ds \right) \right) \right]$$

(2.9)

from the limiting hypothesis on the integral of the function $\gamma_\alpha$, for any arbitrary small real norm-dependent constant $\varepsilon_\alpha \in \mathbf{R}_+$, there exists a finite $t_0 \in \mathbf{Z}_+$ ( only for a simple constructive proof easily extendable to $t_0 \in \mathbf{R}_+$), dependent on $\varepsilon_\alpha$ and the given $\alpha$-norm, such that one gets from (2.9) by taking initial conditions at $t_0$:

$$\|x(t)\|_\alpha \le K_1(\alpha) \frac{(t-t_0)^{\nu-1}}{\nu!} \left[ e^{\mu(t-t_0)} \|x(t_0^+)\|_\alpha + \left|\frac{e^{\mu(t-t_0)}(1-e^{\mu h})}{\mu}\right| (|\phi|_\alpha + K_{0\alpha}) + \beta \varepsilon_\alpha \sup_{t-t_0 \le \tau \le t}(|x_\tau|_\alpha) \right]$$

(2.10)

with $\mu$ being the real part of a characteristic zero of $\hat{T}^{-1}(s)$ of multiplicity $\nu$ and $j_t = max(z \in \mathbf{Z}_{0+} : t \ge j_t)$ is dependent on t. Note that if the solution x(t) is unbounded for the given initial conditions then there exist, by construction, a finite and $t_0 \in \mathbf{Z}_+$ and a subsequence $x(t_k)$ valued at the real increasing sequence $\{t_k\}_0^\infty$ (then $t_k \to \infty$ as $k \to \infty$) such that $\|x(t_k)\|_\alpha = \sup_{t-t_0 \le \tau \le t_k}(|x_\tau|_\alpha)$ so that from (2.10) and for some bounded vector function $g_\alpha \in PC^{(0)}(\mathbf{R}_{0+}, \mathbf{R}_{0+})$:

$$\|x(t_k)\|_\alpha = \sup_{t_k-t_0 \le t \le t_k}(|x_t|_\alpha) \le (1-\beta\varepsilon_\alpha K_1(\alpha))^{-1} \left( K_1(\alpha) \frac{(t_k-t_0)^{\nu-1}}{\nu!} \left[ e^{\mu(t_k-t_0)} \|x(t_0^+)\|_\alpha + \left|\frac{e^{\mu(t_k-t_0)}(1-e^{\mu h})}{\mu}\right| (|\phi|_\alpha + K_{0\alpha}) \right] \right)$$

(2.11)

provided that $\varepsilon_\alpha$ is sufficiently small to guarantee $1 > \varepsilon_\alpha K_1(\alpha)$ in the case that $\beta = 1$ and independently of $\varepsilon_\alpha$ if $\beta = 0$. Furthermore, if $\mu \ne 0$ then $\mu > 0$ if the solution is unbounded since, otherwise, $\sup_{t_k-t_0 \le t \le t_k}(|x_t|_\alpha)$ is bounded from (2.7) which contradicts the made assumption that it is unbounded. The equivalent contrapositive proposition to the last above one is that if $\sup_{t_k-t_0 \le t \le t_k}(|x_t|_\alpha)$ is uniformly bounded then $\mu \le 0$. Equivalently, if furthermore $\mu = 0$, then $\nu = 1$ (i.e. $\mu$ is the real part of a simple real characteristic zero of $\hat{T}^{-1}(s)$ associate with the limiting equation (2.2) or there are two simple complex conjugate ones with real part $\mu$). Otherwise, some unbounded lower-bound may be obtained similarly to (2.10) with the replacement of one of the plus signs in the right-hand-side terms to a minus sign affecting some unbounded term caused by $\frac{(t-t_0)^{\nu-1}}{\nu!} \to \infty$ as $t \to \infty$ if $\nu \ne 1$. This implies that the solution is unbounded which contradicts the fact that it is bounded. Note from (2.10) that if $\mu > 0$ then real increasing sequence $\{t_k\}_0^\infty$ of (2.10):

$$\sup_{t_j-t_0 \le t \le t_k}(|x_{t_j}|_\alpha) = K_1(\alpha) \frac{(t_k-t_0)^{\nu-1}}{\nu!} \left[ e^{\mu(t_k-t_0)} \|x(t_0^+)\|_\alpha + \left|\frac{e^{\mu(t_k-t_0)}(1-e^{\mu h})}{\mu}\right| (|\phi|_\alpha + K_{0\alpha}) + \beta \varepsilon_\alpha \sup_{t_k-t_0 \le \tau \le t_k}(|x_{t_k}|_\alpha) \right] - g_\alpha(t_0, t_k)$$

$$\le (1-\beta\varepsilon_\alpha K_1(\alpha))^{-1} \left( K_1(\alpha) \frac{(t_k-t_0)^{\nu-1}}{\nu!} \left[ e^{\mu(t_k-t_0)} \|x(t_0^+)\|_\alpha + \left|\frac{e^{\mu(t_k-t_0)}(1-e^{\mu h})}{\mu}\right| (|\phi|_\alpha + K_{0\alpha}) \right] - g_\alpha(t_0, t_k) \right)$$

(2.12)

which takes the form $\sup_{t_k - t_0 \leq t \leq t_k}(|x_t|_\alpha) = \frac{(t_k - t_0)^{\nu-1}}{\nu!} e^{\mu(t_k - t_0)} M - g_\alpha(t_0, t_k)$, where $M \in \mathbf{R}_{0+}$ depends on $t_0$ (finite), $K_{1\alpha}$, $K_{0\alpha}$, $\beta$, $\varepsilon_\alpha$, $|\phi_\alpha|$ and the $\alpha$-norm, for some bounded vector function $g_\alpha \in PC^{(0)}(\mathbf{R}_{0+}, \mathbf{R}_{0+})$ which depends on $t_0$, the initial conditions, $\beta$ and $\varepsilon_\alpha$ provided that $\varepsilon_\alpha$ is sufficiently small to guarantee $1 > \varepsilon_\alpha K_1(\alpha)$ in the case that $\beta = 1$ and independently of $\varepsilon_\alpha$ if $\beta = 0$. Assume that the solution x(t) is not a trivial solution what is guaranteed if $\lim_{t \to \infty} e^{bt} x(t) \neq 0$, $\forall b \in \mathbf{R}$. Then, it follows from (2.12) that $\limsup_{t_k - t_0 \leq t \leq t_k}\left(\frac{\ln|x_t|_\alpha}{t^\nu}\right) = \mu > 0$ for sufficiently large $t_0 \in \mathbf{Z}_{0+}$, irrespective of the $\alpha$-norm, since $\lim_{t \to \infty} \frac{\ln g_\alpha(t_0, t)}{t^\nu} = 0$ for any $\alpha$-norm. By taking $\mathbf{Z}_{0+} \ni t_0 \to \infty$, it follows that $\lim_{t \to \infty} \frac{\ln|x_t|_\alpha}{t^\nu} = \mu > 0$. The result may also be extended to the case $\mu = 0$ since then either the solution is unbounded for some initial conditions and multiplicity $\nu > 1$ of the characteristic zero of $\hat{T}^{-1}(s)$ whose real part is $\mu$, or it is bounded (in particular, always if $\nu = 1$). As a result, if $\mu \geq 0$ and there is no $b \in \mathbf{R}$ such that $e^{bt} x(t)$ converges to zero as $t \to \infty$ then $\lim_{t \to \infty} \frac{\ln|x_t|_\alpha}{t^\nu} = \mu$ and $\lim_{t \to \infty} \frac{\ln|x_t|_\alpha}{t^{\nu+\ell}} = \mu$; $\forall \ell \in \mathbf{Z}_+$ since $\int_t^{t+1} \frac{t^{\vartheta_k - 1}}{\vartheta_k!} e^{\beta \sigma_k t} \gamma_\alpha(s) ds \to 0$ as $t \to \infty$ for all the characteristic zeros of the limiting equation (2.2). If (2.7) holds, in particular, with $K_{0\alpha} = 0$ then the above result is also valid from (2.12) for a negative value of $\mu$. Property (ii) has been proved.

(iiii) If (2.7) does not hold for $K_{0\alpha} = 0$ and all the characteristic zeros of the limiting equation (2.2) have negative real parts then it follows by using close reasoning to that used in (ii) that the solution cannot converge asymptotically to zero but it is uniformly bounded from (2.12) since $\int_t^{t+1} \frac{t^{\vartheta_k - 1}}{\vartheta_k!} e^{\beta \sigma_k t} \gamma_\alpha(s) ds \to 0$ as $t \to \infty$ and such an integral is bounded, $\forall t \in \mathbf{R}_+$. Thus, $\mu_k = 0$, $\forall k \geq k_1$ and $\mu_{k_1}$ is not the real part of a characteristic zero of the limiting equation (2.2) since it is not a negative real number. □

The real limit $\mu_{k_1}$ of Theorem 2.1 (ii)-(ii), provided that it exists, is called the *strict Lyapunov exponent* of the solution of (2.2)-(2.3) with the perturbation function $f(t, x_t)$, subject to the hypotheses of Theorem 2.1, which is the real part of an eigenvalue (or characteristic zero) of the limiting equation (2.1) if either it is positive or if it takes any arbitrary value in the case that Eq. (2.7) holds for $K_{0\alpha}=0$ [Theorem 2.1 (ii)]. If all the characteristic zeros of (2.1) have negative real parts but (2.7) is not fulfilled with $K_{0\alpha}=0$ then the strict Lyapunov exponent, if it exists, is zero so that it is not the real part of a characteristic zero of the limiting equation (2.1) [Theorem 2.1(iii)]. The main extension of Theorem 2.1 for the very general functional differential equation (2.2)-(2.3) with respect to parallel previous results (see [1, Chapter IV, Theorem 5 for ordinary differential equations; 2, Theorem 1.1, for functional differential equations] and [19]) is that the perturbation function in (2.7) is not vanishing for bounded solutions or slightly growing solutions since any bounded functions are primarily admitted as perturbations in (2.2). The extension concerning the result in [1] is restricted to the form of (2.1) which involves a wide type of delayed dynamics involving any finite numbers of point delays, finite distributed delays and delays generated by Volterra- type dynamics.

A notation for the subsequent lemma and theorem is the following (see [18, Chapter 7]). If $\Lambda$ is a finite set of eigenvalues of (2.1), then $P_\Lambda$ and $Q_\Lambda$ denote the generalized eigenspace associated with $\Lambda$ and the corresponding complementary subspace of $C(\mathbf{R}_{0+})$, respectively. The phase space $C(\mathbf{R}_{0+})$ is decomposed by $\Lambda$ into the direct sum $C(\mathbf{R}_{0+}) = P_\Lambda(\mathbf{R}_{0+}) \oplus Q_\Lambda(\mathbf{R}_{0+})$. The projections of the solution $x \in C(\mathbf{R}_{0+})$ of (2.2), subject to (2.3), for any initial condition $\phi \in C_e(-h)$, onto the above subspaces are denoted by $x^{P_\Lambda(\mathbf{R}_{0+})}$ and $x^{Q_\Lambda(\mathbf{R}_{0+})}$, respectively, $\forall t \in \mathbf{R}_{0+}$. Note that although the initial conditions of (2.2)-(2.3) are in general in $C_e(-h)$, the corresponding unique solution of (2.2), subject to (2.3), for $t \in \mathbf{R}_{0+}$ are in $C(\mathbf{R}_{0+})$. The whole solutions in $\mathbf{R}_{-h}$ which includes any given initial condition $\phi \in C_e(-h)$ then satisfying $x(t) = \phi(t), \forall t \in [-h, 0]$, and the differential equation (2.2), subject to (2.3), for $t \in \mathbf{R}_{0+}$ are in $C_e(\mathbf{R}_{-h}) = P_\Lambda(\mathbf{R}_{-h}) \oplus Q_{\Lambda e}(\mathbf{R}_{-h})$ where $Q_{\Lambda e}(\mathbf{R}_{-h})$ is the complementary subspace of $P_\Lambda(\mathbf{R}_{-h})$ in $C_e(\mathbf{R}_{-h})$. The projections of the solution onto those subspaces are $x^{P_\Lambda(\mathbf{R}_{-h})}$ and $x^{Q_{\Lambda e}(\mathbf{R}_{-h})}$, respectively, $\forall t \in \mathbf{R}_{-h}$. The following technical result is direct without proof.

**Lemma 2.2.** Assume that the initial condition of (2.2)-(2.3) is $x(t) = \phi(t), \forall t \in [-h, 0]$ for any given $\phi \in C_e(-h)$. The unique solution of (2.2), subject to (2.3) on $\mathbf{R}_{0+}$, and identified with $\phi(t)$ $\forall t \in [-h, 0]$, satisfies with unique decompositions:

$$x_{t+h} = x_{t+h}^{P_\Lambda(\mathbf{R}_{0+})} + x_{t+h}^{Q_\Lambda(\mathbf{R}_{0+})}; \ \forall t \in \mathbf{R}_{0+}$$

$$x_t = x_t^{P_\Lambda(\mathbf{R}_{-h})} + x_t^{Q_{\Lambda e}(\mathbf{R}_{-h})}; \ \forall t \in \mathbf{R}_{0+}$$

$$x_0 = \phi_0^{P_\Lambda(\mathbf{R}_{-h})} = \phi_0^{P_\Lambda(\mathbf{R}_{0+})} = \phi_0^{P_\Lambda(\mathbf{R}_{0+})} = x_0^{P_\Lambda(\mathbf{R}_{-h})} + x_0^{Q_{\Lambda e}(\mathbf{R}_{-h})} = \phi_0^{P_\Lambda(\mathbf{R}_{-h})} + \phi_0^{Q_{\Lambda e}(\mathbf{R}_{-h})} \qquad \square$$

The meaning of Lemma 2.2 is that for $t \geq 0$, x(t) is decomposed uniquely as a sum of a function in $P_\Lambda(\mathbf{R}_{0+})$ and another one in its complementary in $C(\mathbf{R}_{0+})$, even for initial conditions in $C_e(-h)$, rather than in the more restrictive set $C(-h)$. However, for $t \leq 0$, since $x_0 = \phi_0$ for any given $\phi \in C_e(-h)$, the complementary set $Q_{\Lambda e}(\mathbf{R}_{-h})$ of $P_\Lambda(\mathbf{R}_{-h})$ in $C(\mathbf{R}_{-h})$ replaces $Q_\Lambda(\mathbf{R}_{0+})$. Note that $x_t = \phi_t^{P_\Lambda(\mathbf{R}_{-h})} + \phi_t^{Q_{\Lambda e}(\mathbf{R}_{-h})}, \forall t \in \mathbf{R}_{0+}$ is untrue except for $\phi \in C(-h)$.

**Theorem 2.3**. Let x be a solution of (1.2)-(1.3) satisfying the hypotheses of Theorem 2.1 with a finite strict Lyapunov exponent $\mu(x) = \mu_{k_1} = \mu$. Consider generalized eigenspaces $P_0(\boldsymbol{R}_{0+}) = P_{\Lambda_0}(\boldsymbol{R}_{0+})$, $P_1(\boldsymbol{R}_{0+}) = P_{\Lambda_1}(\boldsymbol{R}_{0+})$ and $Q(\boldsymbol{R}_{0+}) = Q_\Lambda(\boldsymbol{R}_{0+})$ for $t \in \boldsymbol{R}_{0+}$ and, also, generalized eigenspaces $P_0(\boldsymbol{R}_{-h}) = P_{\Lambda_0}(\boldsymbol{R}_{-h})$, $P_1(\boldsymbol{R}_{-h}) = P_{\Lambda_1}(\boldsymbol{R}_{-h})$ and $Q_e(\boldsymbol{R}_{-h}) = Q_{\Lambda_e}(\boldsymbol{R}_{-h})$ for $t \in \boldsymbol{R}_{+h}$, where the spectral sets $\Lambda_0$, $\Lambda_1$ and $\Lambda$ each generating the two corresponding eigenspaces, are defined by:

$$\Lambda_0 = \Lambda_0(\mu) := \{\lambda : \det \hat{T}^{-1}(\lambda) = 0, \text{ Re}\,\lambda = \mu\} \tag{2.13}$$

$$\Lambda_1 = \Lambda_1(\mu) := \{\lambda : \det \hat{T}^{-1}(\lambda) = 0, \text{ Re}\,\lambda > \mu\} \tag{2.14}$$

$$\Lambda = \Lambda_0 \cup \Lambda_1 = \Lambda(\mu) := \{\lambda : \det \hat{T}^{-1}(\lambda) = 0, \text{ Re}\,\lambda \geq \mu\} \tag{2.15}$$

Then, the following properties hold:

(i)

(i.1) $\Lambda_0 \neq \Phi$ and $\Lambda \neq \Phi$ if $\mu \in \boldsymbol{R}_{0+}$

(i.2) $\Lambda_0 \neq \Phi$ if $\mu \in \boldsymbol{R}$ and, furthermore, Eq. (2.7) holds with $K_{0\alpha} = 0$.

(i.3) $\Lambda = \Lambda_0 = \Lambda_1 = \Phi$ if all the eigenvalues of (2.1) have negative real parts and, furthermore, Eq.(2.7) does not hold with $K_{0\alpha} = 0$.

(i.4) $\Lambda_1 = \Phi$ if any of the following conditions hold:

   (1) No eigenvalue of (2.1) is in $\boldsymbol{R}_{0+}$.

   (2) No eigenvalue of (2.1) is in $\boldsymbol{R}_+$ and, furthermore, (2.7) does not hold with $K_{0\alpha} = 0$.

   (3) Eq. (2.7) holds with $K_{0\alpha} = 0$.

(ii) The solution of (2.2) under arbitrary initial conditions $\phi \in C_e(-h)$, subject to a perturbation function (2.3), satisfies:

$$x_{t+h} = x_{t+h}^{P_0(\boldsymbol{R}_{0+})} + x_{t+h}^{P_1(\boldsymbol{R}_{0+})} + x_{t+h}^{Q(\boldsymbol{R}_{0+})}; \quad \forall t \in \boldsymbol{R}_{0+} \tag{2.16}$$

$$x_t = x_t^{P_0(\boldsymbol{R}_{-h})} + x_t^{P_1(\boldsymbol{R}_{-h})} + x_t^{Q_e(\boldsymbol{R}_{-h})}; \quad \forall t \in \boldsymbol{R}_{0+} \tag{2.17}$$

$$x_0 = \phi_0^{P_0(\boldsymbol{R}_{-h})} = x_0^{P_0(\boldsymbol{R}_{0+})} = \phi_0^{P_0(\boldsymbol{R}_{0+})} = x_0^{P_0(\boldsymbol{R}_{-h})} + x_0^{P_1(\boldsymbol{R}_{-h})} + x_0^{Q_e(\boldsymbol{R}_{-h})} = \phi_0^{P_0(\boldsymbol{R}_{-h})} + x_0^{P_1(\boldsymbol{R}_{-h})} + \phi_0^{Q_e(\boldsymbol{R}_{-h})}$$

$$\tag{2.18}$$

(iii) The solution of (2.2) under arbitrary initial conditions $\phi \in C_e(-h)$, subject to a perturbation function (2.3), satisfies:

$$x_{t+h}^{P_1(\boldsymbol{R}_{0+})} = O\left(\left|x_{t+h}^{P_0(\boldsymbol{R}_{0+})}\right|\right),\ x_{t+h}^{Q(\boldsymbol{R}_{0+})} = O\left(\left|x_{t+h}^{P_0(\boldsymbol{R}_{0+})}\right|\right),\ x_{t+h} = O\left(\left|x_{t+h}^{P_0(\boldsymbol{R}_{0+})}\right|\right);\ \forall t \in \boldsymbol{R}_{0+}$$

$$\tag{2.19}$$

$$x_t^{P_1(\boldsymbol{R}_{-h})} = O\left(\left|x_t^{P_0(\boldsymbol{R}_{-h})}\right|\right),\ x_t^{Q_e(\boldsymbol{R}_{-h})} = O\left(\left|x_t^{P_0(\boldsymbol{R}_{-h})}\right|\right),\ x_t = O\left(\left|x_t^{P_0(\boldsymbol{R}_{-h})}\right|\right);\ \forall t \in \boldsymbol{R}_{0+}$$

$$\tag{2.20}$$

Furthermore, if Eq. (2.7) holds with $K_{0\alpha} = 0$ then, as $t \to \infty$:

$$x_{t+h}^{P_1(\boldsymbol{R}_{0+})} = o\left(\left|x_{t+h}^{P_0(\boldsymbol{R}_{0+})}\right|\right),\ x_{t+h}^{Q(\boldsymbol{R}_{0+})} = o\left(\left|x_{t+h}^{P_0(\boldsymbol{R}_{0+})}\right|\right),\ x_{t+h} = o\left(\left|x_{t+h}^{P_0(\boldsymbol{R}_{0+})}\right|\right);\ \forall t \in \boldsymbol{R}_{0+}$$

$$\tag{2.21}$$

$$x_t^{P_1}(\mathbf{R}_{-h}) = o\left(\left|x_t^{P_0}(\mathbf{R}_{-h})\right|\right), \; x_t^{Q_e}(\mathbf{R}_{-h}) = o\left(\left|x_t^{P_0}(\mathbf{R}_{-h})\right|\right), \; x_t = o\left(\left|x_t^{P_0}(\mathbf{R}_{-h})\right|\right); \; \forall t \in \mathbf{R}_{0+}$$

(2.22)

(iv) The solution of the limiting equation (2.1) satisfies (2.21)-(2.22) as $t \to \infty$.

(v) The solution of (2.2), under arbitrary initial conditions $\phi \in C_e(-h)$ and subject to a perturbation function (2.3), satisfies:

$$x_{t+h} = x_{t+h}^{P_0} + x_{t+h}^{P_1} + x_{t+h}^{Q}; \; \forall t \in \mathbf{R}_{0+} \tag{2.23}$$

$$x_t = x_t^{P_0} + x_t^{P_1} + x_t^{Q_e}; \; \forall t \in \mathbf{R}_{0+} \tag{2.24}$$

which is identical to

$$x_t = x_t^{P_0} + x_t^{P_1} + x_t^{Q}; \; \forall t \in R_{-h} \tag{2.25}$$

under the restriction $\phi \in C(-h)$ for the initial conditions with $x(t) = \phi(t)$. Also,

$$x_t^{P_1} = O\left(\left|x_t^{P_0}\right|\right), \; x_t^{Q} = O\left(\left|x_t^{P_0}\right|\right), \; x_t = O\left(\left|x_t^{P_0}\right|\right); \; \forall t \in \mathbf{R}_{0+} \tag{2.26}$$

$\forall \phi \in C_e(-h)$. If Eq. (2.7) holds with $K_{0\alpha} = 0$ then, as $t \to \infty$:

$$x_{t+h}^{P_1} = o\left(\left|x_{t+h}^{P_0}\right|\right), \; x_{t+h}^{Q} = o\left(\left|x_{t+h}^{P_0}\right|\right), \; x_{t+h} = o\left(\left|x_{t+h}^{P_0}\right|\right); \; \forall t \in \mathbf{R}_{0+} \tag{2.27}$$

which leads to

$$x_t^{P_1} = o\left(\left|x_t^{P_0}\right|\right), \; x_t^{Q} = o\left(\left|x_t^{P_0}\right|\right), \; x_t = o\left(\left|x_t^{P_0}\right|\right); \; \forall t \in \mathbf{R}_{0+} \tag{2.28}$$

under the restriction $\phi \in C(-h)$ for the initial conditions with $x(t) = \phi(t)$. The solution of the limiting equation (2.1) satisfies (2.27)-(2.28) as $t \to \infty$.

*Proof*: Property (i) is a direct consequence of Theorem 2.1 [(ii)-(iii)]. Property (ii) is a direct consequence of Lemma 2.2 since $\Lambda_0$ and $\Lambda_1$ are disjoint sets what implies that $C(\mathbf{R}_{0+}) = P_0(\mathbf{R}_{0+}) \oplus P_1(\mathbf{R}_{0+}) \oplus Q_\Lambda(\mathbf{R}_{0+})$ and $C_e(\mathbf{R}_{-h}) = P_0(\mathbf{R}_{-h}) \oplus P_1(\mathbf{R}_{-h}) \oplus Q_{\Lambda e}(\mathbf{R}_{-h})$. Eqs. (2.19)- (2.20) are a direct consequence of Property (ii). Eqs. (2.21)- (2.22) are a direct consequence of (2.19)-(2.20) if (2.7) holds for $K_{0\alpha} = 0$ so that Property (iii) follows. Property (iv) follows from Property (iii) as particular case for $f(t, x_t) = 0$; $\forall t \in \mathbf{R}_{-h}$ in Eq. (2.3). Property (v) is a direct consequence of Properties (i)–(iv) In particular, the relative growing "O"-properties of the various parts of the solution of (2.2)-(2.3) are embedded from Property (iii) into similar properties for the solution strings of length h. The part of Property (v) concerning the relative growing "o"-properties of the various parts of the solution of (2.2)-(2.3) and that concerning the limiting equation follows under a close reasoning. □

Note that in Theorem 2.3, the various results obtained for "Landau´s small -o" notation, referred to limits as $t\to\infty$ imply, as usual, that parallel results for "Landau´s big -O" notation stand for all $t\in R_{0+}$ but the converse is not true. The results concerning "Landau´s big -O" notation in Theorem 2.3 (iii) for the perturbed functional equation (2.2)-(2.3) are new for the studied class of functional equations, related to the background literature, since the perturbation function is permitted to take bounded nonzero values even if the limiting equation is globally asymptotically stable and it is not requested to grow asymptotically at most linearly with $x_t$. The results concerning "Landau´s big -O" notation imply that the solution of the perturbed functional equation is uniformly bounded for any bounded function of initial conditions of the given class for all time so that the functional differential equation is globally uniformly Lyapunov´s stable provided that the perturbation (2.3) satisfies the given hypotheses. A technical result concerning the boundedness of the evolution operator, which will be then useful to derive further results, and stability properties of the differential systems (2.1) and (2.2)-(2.3) follows:

**Theorem 2.4**. The following properties hold:

(i) The evolution operator of the limiting functional differential equation (2.1) satisfies the subsequent relations:

$$\|T(t,0)\|_\alpha \le K_1(\alpha)\|I_n\|_\alpha \max\left(1, \frac{t^{\nu-1}}{\nu!} e^{\mu t}\right); \forall t \in R_{0+} \tag{2.29}$$

$$|T(t,0)|_\alpha \le K_1(\alpha)\|I_n\|_\alpha \max\left(1, \frac{t^{\nu-1}}{\nu!} e^{\mu t} \max\left(1, \sum_{i=0}^{\nu-2} t^{1-\nu+i} h^{\nu-1-i} e^{\mu h}\right)\right); \forall t \in R_{0+} \tag{2.30}$$

$$\|\dot{T}(t,0)\|_\alpha \le K_2(\alpha) \sup_{0\le\tau\le t}\left|T_\tau^{C(R_{0+})}\right|_\alpha; \forall t \in R_{0+} \tag{2.31}$$

$$|T(t,0)|_\alpha \le K_2(\alpha) \sup_{0\le\tau\le t}\left|T_\tau^{C(R_{0+})}\right|_\alpha \tag{2.32}$$

$$\le K_1(\alpha)K_2(\alpha)\|I_n\|_\alpha \max\left(1, \frac{t^{\nu-1}}{\nu!} e^{\mu t} \max\left(1, \sum_{i=0}^{\nu-2} t^{1-\nu+i} h^{\nu-1-i} e^{\mu h}\right)\right); \forall t \in R_{0+} \tag{2.33}$$

*Proof*: (i) The evolution operator satisfies the limiting functional differential equation (2.1):

$$\dot{T}(t,0) = \sum_{i=0}^{m} A_i T(t-h_i,0) + \sum_{i=0}^{m'} \int_0^t d\alpha_i(\tau) A_{\alpha_i} T(t-h'_i,\tau) + \sum_{i=m'+1}^{m'+m''} \int_{h'_i}^0 d\alpha_i(\tau) A_{\alpha_i} T(t,\tau) \tag{2.34}$$

for $t\in R_{0+}$ subject to initial conditions $T(0,0)=I_n$ (i.e. the n-th identity matrix) and $T(t,0)=0$, $t\in[-h,0)$. Thus, it satisfies also the unforced (2.9) (i.e. for $\gamma_\alpha = |\phi|_\alpha = K_{0\alpha} = 0$). This leads directly to (2.29). Eq. (2.30) follows by using the Newton binomial to expand $\frac{(t+h)^{\nu-1}}{\nu!} e^{\mu(t+h)}$ and the fact that the maximum of the real exponential function within the real interval $[0,t]$ is reached at the boundary. Eq. (2.31) follows by inspection from (2.34) for some norm-dependent $K_2(\alpha) \in \mathbf{R}_+$ which depends on the various matrices of parameters of the limiting functional differential equation (2.1). Eq. (2.32) follows from Eqs. (2.31) and (2.34). Finally, Eq. (2.33) follows from Eqs. (2.32) and (2.30). property (i) has been proved.

(ii) For sufficiently small constant $(\beta \varepsilon_\alpha)$, the evolution operator as a function of time is of exponential order whose norm time-function satisfies:

$$\|T(t,0)\|_\alpha = \sup_{t_k - t_0 \leq t \leq t_k} (|x_t|_\alpha) \leq (1 - \beta \varepsilon_\alpha K_1(\alpha))^{-1} \left( K_1(\alpha) \frac{(t_k - t_0)^{v-1}}{v!} \left[ e^{\mu(t_k - t_0)} \|I_n\|_\alpha \right] \right)$$

(2.35)

which converges exponentially to zero as $t \to \infty$ if the strict Lyapunov exponent $\mu$ is negative. In this case, the limiting differential functional equation is globally uniformly exponentially Lyapunov´s stable whose solution satisfies asymptotically:

$$\|x(t_k)\|_\alpha = \sup_{t_k - t_0 \leq t \leq t_k} (|x_t|_\alpha) \leq (1 - \beta \varepsilon_\alpha K_1(\alpha))^{-1} \left( K_1(\alpha) \frac{(t_k - t_0)^{v-1}}{v!} \left[ e^{\mu(t_k - t_0)} \|x(t_0^+)\|_\alpha + \left| \frac{e^{\mu(t_k - t_0)}(1 - e^{\mu h})}{\mu} \right| (|\phi|_\alpha) \right] \right)$$

(2.36)

so that it converges exponentially to zero as $t \to \infty$ for any admissible function of initial conditions. The differential equation (2.2), subject to (2.3) is globally uniformly Lyapunov´s stable if $\mu \leq 0$ and its solution satisfies for large t:

$$\|x(t_k)\|_\alpha = \sup_{t_k - t_0 \leq t \leq t_k} (|x_t|_\alpha) \leq (1 - \beta \varepsilon_\alpha K_1(\alpha))^{-1} \left( K_1(\alpha) \frac{(t_k - t_0)^{v-1}}{v!} \left[ e^{\mu(t_k - t_0)} \|x(t_0^+)\|_\alpha + \left| \frac{e^{\mu(t_k - t_0)}(1 - e^{\mu h})}{\mu} \right| (|\phi|_\alpha + K_{0\alpha}) \right] \right)$$

(2.37)

and converges exponentially to zero (i.e. it is globally uniformly exponentially Lyapunov´s stable) if $\mu < 0$ and the perturbation function has an upper-bounding function with $K_{0\alpha} = 0$.

(ii) It follows directly from (2.9). □

The evolution operator $T: \mathbf{R}_{0+} \times \mathbf{C}^n \to \mathbf{C}^n$ which makes explicit the solutions of the limiting equation Eq. (2.4) and the perturbed one Eq. (2.5) for each function of initial conditions. Then, let $(T_s(t))_{t \in \mathbf{R}_{0+}}$ be the solution semigroup of the linear autonomous equation (2.1), which is unique for $t \in \mathbf{R}_{0+}$ for each $\phi \in C_e(-h)$ and whose infinitesimal generator is $A$ satisfying $\dot{\varphi} = A\varphi$, $\forall \varphi \in Dom(A) := \{\varphi \in C_e(\mathbf{R}_{-h}) : \exists \dot{\varphi} \in C(\mathbf{R}_{0+}), \forall t \in \mathbf{R}_{0+} \wedge \dot{\varphi}(0) = L\varphi(0)\}$. Thus, the string $x_t(\phi) = (T_s \phi)(t)$ of the solution of the limiting functional differential equation (2.1) within $[t-h, t]$ is defined from (2.4) as follows:

$$x_t(\phi) = (T_s \phi)(t) := T(t-\theta, 0)x(0^+) + \int_{-h}^{0} T(t-\theta, \tau)\phi(\tau)U(\tau)d\tau, \forall \theta \in [0, \min(t, h)]; \text{ and}$$

$$x_t(\phi) = (T_s \phi)(t) = 0, \forall \theta \in ]0, \min(t, h)[ := \overline{[0, \min(t, h)]} \cap \mathbf{R}_{0+}$$ (2.38)

$\forall t \in \mathbf{R}_{0+}$ ; and the corresponding solution string of the perturbed functional differential equation (2.2) – (2.3) is then defined follows:

$$x_t(\phi) = (T_s \phi)(t) := T(t-\theta, 0)x(0^+) + \int_{-h}^{0} T(t-\theta, \tau)\phi(\tau)U(\tau)d\tau + \int_0^t T(t-\theta, \tau)f(\tau, x_\tau)d\tau$$

, $\forall \theta \in [0, min(t,h)]$; and

$$x_t(\phi) = (T_s \phi)(t) = 0 \ , \ \forall \theta \in \overline{[0, min(t,h)]} \cap \mathbf{R} \ ; \ \forall t \in \mathbf{R}_{0+} \tag{2.39}$$

The transposed equation associated with (2.1) is:

$$\dot{y}(t) = L^* y_t := \sum_{i=0}^{m} y(t-h_i)A_i^* + \sum_{i=0}^{m'} \int_0^t y(t-\tau-h_i')A_{\alpha_i}^* d\alpha_i^*(\tau) + \sum_{i=m'+1}^{m'+m''} \int_{t-h_i'}^{t} y(\tau)A_{\alpha_i}^* d\alpha_i^*(\tau)$$

(2.40)

where the superscript * denotes the adjoint operators of the corresponding un- superscripted ones. In particular, for matrices, it denotes the conjugate transposes of the corresponding un- superscripted ones. Thus, y(t) is a n- dimensional complex row vector. The phase space for (2.40) on $\mathbf{R}_{0+}$ is $C'(\mathbf{R}_{0+}) := C(\mathbf{R}_{0+}, \mathbf{C}^{n^*})$. Corresponding spaces of functions taking into account the more general spaces for initial conditions are $C'(h) := C([0,h], \mathbf{C}^{n^*})$, $C_e'(h) := C([0,h], \mathbf{C}^{n^*})$ and $C'(\mathbf{R}_{-h}) := C(\mathbf{R}_{-h}, \mathbf{C}^{n^*})$. Let $\Lambda$ be a finite set of eigenvalues of (2.1) and let $\Phi_\Lambda$ be a basis for the generalized eigenspace $P_\Lambda$, [1-2]. Then, there exists a square n-matrix $B_\Lambda$, with $sp(B_\Lambda) = sp(\Lambda)$, such that the subsequent relations hold:

$$A\Phi_\Lambda = \Phi_\Lambda B_\Lambda, \quad \Phi_\Lambda(\tau) = \Phi_\Lambda(0)e^{B_\Lambda \tau} (\forall \tau \in [-h, 0]), \quad T(t,0)\Phi_\Lambda = \Phi_\Lambda e^{B_\Lambda t} \tag{2.41}$$

The relations (2.41) yield via direct computation Property (i) of the subsequent result since $B_\Lambda$ commutes with $e^{B_\Lambda t}$. Property (ii) is a direct consequence of (2.34) subject to $T(0,0) = I_n$ and $T(t,0) = 0$ for $t < 0$.

**Proposition 2.5**. The two following properties hold:
(i) The following relations hold, $\forall t \in \mathbf{R}_{0+}$ :

$$T(t,0)\Phi_\Lambda B_\Lambda = \Phi_\Lambda e^{B_\Lambda t} B_\Lambda = \Phi_\Lambda B_\Lambda e^{B_\Lambda t} = A\Phi_\Lambda e^{B_\Lambda t} = AT(t,0)\Phi_\Lambda \tag{2.42}$$

$$\dot{T}(t,0)\Phi_\Lambda B_\Lambda = \Phi_\Lambda B_\Lambda^2 e^{B_\Lambda t} = A\Phi_\Lambda B_\Lambda e^{B_\Lambda t} = A^2 \Phi_\Lambda e^{B_\Lambda t} = A\Phi_\Lambda e^{B_\Lambda t} B_\Lambda \tag{2.43}$$

$$= \Phi_\Lambda B_\Lambda e^{B_\Lambda t} B_\Lambda = T(t,0)\Phi_\Lambda B_\Lambda^2 = T(t,0)A\Phi_\Lambda B_\Lambda = A^2 T(t,0)\Phi_\Lambda \tag{2.44}$$

$$= AT(t,0)\Phi_\Lambda B_\Lambda = AT(t,0)A\Phi_\Lambda \tag{2.45}$$

(ii) The evolution operator of the solution of (2.1) is uniquely given by:

$$T(t,0) = e^{A_0 t}\left(I_n + \left[\int_0^t e^{-A_0 \tau}\left(\sum_{i=1}^{m} A_i T(t-h_i,\tau) + \sum_{i=0}^{m'} \int_0^\tau d\alpha_i(\theta) A_{\alpha_i} T(\tau-h_i',\theta) + \sum_{i=m'+1}^{m'+m''} \int_{h_i}^{0} d\alpha_i(\theta) A_{\alpha_i} T(\tau,\theta)\right) d\tau\right]\right)$$

(2.46)

$\forall t \in \mathbf{R}_{0+}$ with $T(0,0) = I_n$ and $T(t,0) = 0$ for $t \in [-h, 0]$. □

Eqns. (2.41)-(2.46) are useful for the asymptotic analysis of comparison of the solutions of (2.2)-(2.3) with that of its limiting equation obtained from (2.1) which follows. The solutions of $P_\Lambda$ can be extended to $\forall t \in \mathbf{R}$ by $T(t,0)\Phi_\Lambda a = \Phi_\Lambda e^{B_\Lambda t} a$, where a is of dimension compatible with the order of $\Phi_\Lambda$. Let $Q_\Lambda$ be the complementary eigenspace to $P_\Lambda$. Now, use appropriate notations for the corresponding subspaces on $\mathbf{R}_{0+}$ and their extensions to $\mathbf{R}_{-h}$ to consider more general initial conditions (on $C_e(-h)$) for (2.1) and (2.2)-(2.3) than bounded continuous functions in a Banach space leading to the uniquely defined decompositions $C_e(\mathbf{R}_{-h}) = P_\Lambda(\mathbf{R}_{-h}) \oplus Q_{\Lambda e}(\mathbf{R}_{-h})$ and $C(\mathbf{R}_{0+}) = P_\Lambda(\mathbf{R}_{0+}) \oplus Q_\Lambda(\mathbf{R}_{0+})$. Then, given a function of initial conditions $\phi \in C_e(-h)$ the decomposition $\phi_0 = \phi_t|_{t=0} = \phi_0^{P_\Lambda(\mathbf{R}_{-h})} + \phi_0^{Q_{\Lambda e}(\mathbf{R}_{-h})}$ is unique. Also, the unique solution of (2.1) and that of (2.2), subject to (2.3), are uniquely decomposable in $\mathbf{R}_{0+}$ as

$$x_{t+h} = x_{t+h}^{P_\Lambda(\mathbf{R}_{0+})} + x_{t+h}^{Q_\Lambda(\mathbf{R}_{0+})} \tag{2.47}$$

$$x_{t+h}^{P_\Lambda(\mathbf{R}_{0+})} = \Phi_\Lambda(\Psi_\Lambda, x_{t+h}) \in P_\Lambda(\mathbf{R}_{0+}), \quad x_{t+h}^{Q_\Lambda(\mathbf{R}_{0+})} = x_{t+h} - x_{t+h}^{P_\Lambda(\mathbf{R}_{0+})} \in Q_\Lambda(\mathbf{R}_{0+}) \tag{2.48}$$

$\forall t \in \mathbf{R}_{0+}$, via the direct sum of subspaces $C(\mathbf{R}_{0+}) = P_\Lambda(\mathbf{R}_{0+}) \oplus Q_\Lambda(\mathbf{R}_{0+})$. The whole solution in $\mathbf{R}_{-h}$ including initial conditions defined by $x(t) = \phi(t)$ for $t \in [-h, 0]$ is uniquely decomposable in $\mathbf{R}_{-h}$ as

$$x_t = x_t^{P_\Lambda(\mathbf{R}_{-h})} + x_t^{Q_{\Lambda e}(\mathbf{R}_{-h})} \tag{2.49}$$

$$x_t^{P_\Lambda(\mathbf{R}_{-h})} = \Phi_\Lambda(\Psi_\Lambda, x_t) \in P_\Lambda(\mathbf{R}_{-h}), \quad x_t^{Q_{\Lambda e}(\mathbf{R}_{-h})} = x_t - x_t^{P_\Lambda(\mathbf{R}_{-h})} \in Q_{\Lambda e}(\mathbf{R}_{-h}) \tag{2.50}$$

$\forall t \in \mathbf{R}_{0+}$, via the direct sum of subspaces $C_e(\mathbf{R}_{-h}) = P_\Lambda(\mathbf{R}_{-h}) \oplus Q_{\Lambda e}(\mathbf{R}_{-h})$.

## 3. Asymptotic behavior

The string solution (2.5) of (2.2)-(2.3) for $\theta \in [t-h, t]$, pointwise defined by $x(t) = \phi(t)$, $t \in [-h, 0]$, any given $\phi \in C_e(-h)$, and

$$x(t+\theta) = T(t+\theta, 0) x(0^+) + \int_{-h}^{0} T(t+\theta, \tau) \phi(\tau) U(\tau) d\tau + \int_{0}^{t+\theta} T(t+\theta, \tau) f(\tau, x_\tau) d\tau$$

(3.1)

$; \forall \theta \in [-h, 0]$, $t \geq h$, may be expressed equivalently via the solution semigroup of the limiting equation (2.1) as

$$x_t(\phi, \theta) = x_t^*(\phi, \theta) + \int_0^{t+\theta} d[K(t+\theta, \tau)] f(\tau, x_\tau) = (T_s(t,0)\phi)(\theta) + \int_0^t d[K(t+\theta, \tau)] f(\tau, x_\tau)$$
$$= T(t+\theta, 0)\phi(0^+) + \int_{-h}^{0} T(t+\theta, \tau)\phi(\tau) U(\tau) d\tau + \int_0^{t+\theta} T(t+\theta, \tau) f(\tau, x_\tau) d\tau$$

(3.2)

where $x^*_{(.)}: [t-h, t] \times C^n \times [-h, 0] \to C^n$, defined by,

$$x^*_t(\phi, \theta) = (T_s(t,0)\phi)(\theta) := T(t+\theta)\phi(0^+) + \int_{-h}^{0} T(t+\theta, \tau)\phi(\tau)U(\tau)d\tau \qquad (3.3)$$

, $\forall \theta \in [-h, 0]$, $t \geq h$, and $(T_s(0,0)\phi)(\theta) = x_0(\phi, \theta) = \phi(\theta)$, $\forall \theta \in [-h, 0]$ with $x(0) = \phi(0)$ is the unique solution of the limiting equation (2.1), and the kernel $K(t,.): [0, t] \to C^n$ of $T_s(t,0)$, $\forall t \in \mathbf{R}_{0+}$ is defined by

$$K(t, s)(\theta) = \int_0^s X(t+\theta-\tau)d\tau, \quad \forall \theta \in [-h, 0], \forall t \in \mathbf{R}_{0+} \qquad (3.4)$$

where X is the fundamental matrix of (2.1) with initial values $X_0(0) = I_n$ and $X_0(\theta) = 0$, $\forall \theta \in [-h, 0]$. The following technical result holds:

**Lemma 3.1**. The following relations hold:

$$x_t = x_t(\phi, \theta) = x_t^{P_A} + x_t^{Q_A}, \quad \mathbf{R}_+ \ni t \geq h \qquad (3.5)$$

$$x_t^{P_A} = \left(T_s(t,0)\phi^{P_A}\right)(\theta) + \int_0^t \left(T_s(t,\tau)f(\tau, x_\tau)^{P_A}\right)(\theta)d\tau$$

$$= \left(T_s(t,0)\phi^{P_A}\right)(\theta) + \int_0^t \left(T_s(t,\tau)X_0^{P_A}\right)(\theta)f(\tau, x_\tau)d\tau$$

$$= T(t+\theta, 0)\phi^{P_A}(0^+) + \int_{-h}^{0} T(t+\theta, \tau)\phi^{P_A}(\tau)U(\tau)d\tau + \int_0^t T(t+\theta, \tau)X_0^{P_A}f(\tau, x_\tau)d\tau$$

(3.6)

$$x_t^{Q_A} = \left(T_s(t,0)\phi^{Q_A}\right)(\theta) + \int_0^t \left(T_s(t,\tau)f(\tau, x_\tau)^{Q_A}\right)(\theta)d\tau, \quad \mathbf{R}_+ \ni t \geq h$$

$$= \left(T_s(t,0)\phi^{Q_{Ae}}\right)(\theta) + \int_0^t \left(d[K(t,\tau)]^{Q_A}\right)(\theta)f(\tau, x_\tau), \quad \mathbf{R}_+ \ni t \geq h \qquad (3.7)$$

$$x_t = \left(T_s(t,0)\left(\phi^{P_A} + \phi^{Q_{Ae}}\right)\right)(\theta) + \int_0^t \left(T_s(t,\tau)X_0^{P_A} + \left(d[K(t,\tau)]^{Q_A}\right)\right)(\theta)f(\tau, x_\tau)d\tau$$

$$, \mathbf{R}_+ \ni t \geq h \quad (3.8)$$

where

$$X_0^{P_A} = \Phi_A \Psi_A(0), \quad K(t,\tau)^{Q_A} = K(t,\tau) - \Phi_A(\Psi_A, K(t,\tau)) \qquad (3.9)$$

Also, the following relations hold for $\phi \in C_e(-h)$, $\forall \varepsilon \in \mathbf{R}_+$ being sufficiently small:

$$\left|T(t,0)\phi^{P_A}\right| \leq M_1 e^{(\mu-\varepsilon)t}\left|\phi^{P_A}\right|, \quad \left|T(t,0)X_0^{P_A}\right| \leq M_1 e^{(\mu-\varepsilon)t}; \forall t \in \mathbf{R}_{0-} \qquad (3.10a)$$

$$\left|T(t,0)\phi^{Q_{Ae}}\right| \leq M_1 e^{(\mu-\varepsilon)t}\left|\phi^{Q_{Ae}}\right|, \quad \forall t \in \mathbf{R}_{0+} \qquad (3.10b)$$

for some $M_1 = M_1(\varepsilon) \in \mathbf{R}_+$ irrespective of the multiplicity of the eigenvalue of the limiting equation (2.1) whose real part is $\mu$.

*Proof.* Eqs. (3.5)-(3.7) hold for any $\mathbf{R}_{0+} \ni t \geq h$ from Theorem 2.3, the definition of the set $\Lambda$ in Theorem 2.1, Eq. (2.15), Lemma 2.2 and (2.47)-(2.49). To obtain (3.5)-(3.7), the following identities are used:

$$\int_0^{t+\theta}\left(T_s(t,\tau)f(\tau,x_\tau)^{P_\Lambda}\right)(\theta)d\tau = \int_0^t \left(T_s(t,\tau)X_0^{P_\Lambda}\right)(\theta)f(\tau,x_\tau)d\tau = \int_0^t T(t+\theta,\tau)X_0^{P_\Lambda} f(\tau,x_\tau)d\tau$$

$$\int_0^t \left(T_s(t,\tau)f(\tau,x_\tau)^{P_\Lambda}\right)(\theta)d\tau = \int_0^{t+\theta}\left(T_s(t,\tau)f(\tau,x_\tau)^{P_\Lambda}\right)(\theta)d\tau = \int_0^t \left(T_s(t,\tau)X_0^{P_\Lambda}\right)(\theta)f(\tau,x_\tau)d\tau$$

since $f(\tau,x_\tau) = 0$ for $\tau < 0$ and $T(t,\tau) = 0$ for $\tau > t$ and $\theta \in [-h, 0]$. Eq.(3.8) follows directly by substitution of (3.6)-(3.7) into (3.5). The norm relations (3.9) hold directly from (3.5)-(3.7) through (3.8).

$\square$

It turns out that for any $\phi \in C_e(-h)$, the above relations hold also for any $t \in \mathbf{R}_{-h}$ by replacing $Q_\Lambda \to Q_{\Lambda e}$ in (3.5)-(3.7). Eqs. (3.9) also hold for $t \in \mathbf{R}_{-h}$ since $\phi \in C_e(-h)$ is bounded. The second relation in (3.5) may be rewritten as $\left|T(t,0)\phi^{Q_\Lambda}\right| \leq M_1 e^{(\mu-\varepsilon)t}$ for any $\mathbf{R}_{0+} \ni t \geq h$ and extended to any $t \in \mathbf{R}_{0+}$ if $\phi$ is continuous on its definition domain $[-h,0]$. A direct consequence of Lemma 3.1 Eq. (3.8) and Proposition 2.5 is that if $x_t^{P_\Lambda} = \Phi_\Lambda u(t)$, $\forall t \in \mathbf{R}_{0+}$ then u is a solution of the ordinary differential equation $\dot{u}(t) = B_\Lambda u(t) + \Psi_\Lambda(0)f(t,x_t)$, [2], which is given explicitly by:

$$u(t) = e^{B_\Lambda t}\left(u(0) + \int_0^t e^{-B_\Lambda \tau} \Psi_\Lambda(0)f(\tau,x_\tau)d\tau\right)$$
$$= e^{B_\Lambda t} u(0) + \int_0^t e^{B_\Lambda(t-\tau)} \Psi_\Lambda(0)f(\tau,x_\tau)d\tau$$
$$= \Phi_\Lambda^{-1} T(t,0)\Phi_\Lambda u(0) + \int_0^t \Phi_\Lambda^{-1} T(t,\tau)\Phi_\Lambda \Psi_\Lambda(0)f(\tau,x_\tau)d\tau$$
$$= \Phi_\Lambda^{-1}\left(T(t,0)\Phi_\Lambda u(0) + \int_0^t T(t,\tau)\Phi_\Lambda \Psi_\Lambda(0)f(\tau,x_\tau)d\tau\right)$$

for $u(0) = \int_{-h}^0 e^{B_\Lambda(t-\tau)}\phi(\tau)d\tau$ so that

$$x_t^{P_\Lambda} = \Phi_\Lambda u(t) = T(t,0)\Phi_\Lambda u(0) + \int_0^t T(t,\tau)\Phi_\Lambda \Psi_\Lambda(0)f(\tau,x_\tau)d\tau$$
$$\dot{u}(t) = B_\Lambda \Phi_\Lambda^{-1} T(t,0)\Phi_\Lambda u(0) + \int_0^t B_\Lambda \Phi_\Lambda^{-1} T(t,\tau)\Phi_\Lambda \Psi_\Lambda(0)f(\tau,x_\tau)d\tau + \Psi_\Lambda(0)f(t,x_t)$$
$$= \Phi_\Lambda^{-1}\left(\dot{T}(t,0)\Phi_\Lambda u(0) + \int_0^t \dot{T}(t,\tau)\Phi_\Lambda \Psi_\Lambda(0)f(\tau,x_\tau)d\tau\right) + \Psi_\Lambda(0)f(t,x_t)$$

It is now proved that the asymptotic difference function between some eigensolution of (2.1), i.e. a finite sum of solutions of (2.1) corresponding to the set $\Lambda_0(\mu)$ of the form $p(t)e^{\lambda t}$, where $p$ is a $\mathbf{C}^n$-valued polynomial and $\lambda \in \Lambda_0(\mu)$, and some corresponding solution of (2.2)-(2.3) grows non faster than linearly with the norm of the solution of the limiting equation (2.1), [16]. If $\gamma(t)$ converges to zero exponentially then the asymptotic difference function between both solutions has strict Lyapunov exponent smaller than that of the corresponding limiting eigensolution. As a result, a solution of (2.2)-(2.3) is of the same exponential order as that of its limiting equation for sufficiently large time. Define the set of distinct eigenvalues of the limiting equation (2.1) as :

$$CE := \left\{\lambda_{j\mu_i} \in \mathbf{C} : \operatorname{Re}\lambda_{j\mu_i} = \mu_i, \det\left(\hat{T}^{-1}(\lambda_{\mu_i})\right) = 0 ; \forall j \in In_{CEi} \subset \mathbf{Z}_+, \forall i \in In_{CER} \subset \mathbf{Z}_+\right\}$$

(3.11)

It is obvious that the cardinal of CE may be infinity (since (2.1) is a functional differential equation), but always numerable. The also denumerable subsets $In_{CEi}$ and $In_{CER}$ of the set of positive integers are indicator sets for all distinct members of CE with real part $\mu_i$ and for all the distinct real parts of the members of CE, respectively. It follows that $Card(CE) = \sum_{i \in In_{CER}} Card\, In_{CEi} \leq \chi_0$. The total number of eigenvalues taking account for their multiplicities $\nu_{ij}$; $\forall j \in In_{CEi}$, $\forall i \in In_{CER}$ is

$$\vartheta = \sum_{i \in In_{CER}} \sum_{j \in In_{CEi}} \nu_{ij}\, Card\, In_{CEi} \geq Card(CE) = \sum_{i \in In_{CER}} Card\, In_{CEi} \qquad (3.12)$$

In particular, $Card(CE) = \chi_0$ ($\chi_0$ standing for infinity denumerable cardinal as opposed to a non-numerable infinity cardinal typically denoted by $\infty$) or $Card(CE)$ is finite. The above definition relies on the fact distinct (non real) eigenvalues $\lambda_{j\mu_i}$ of (2.1), $\forall j \in In_{CEi}$ can have identical real part. Similar sets of eigenvalues of (2.1) as those in (2.13)-(2.15) may be defined being associated to each member of CE as follows:

$$\Lambda_0(\mu_i) := \{\lambda \in CE : Re\,\lambda = \mu_i\} \;;\; \Lambda_1(\mu_i) := \{\lambda \in CE : Re\,\lambda > \mu_i\} \qquad (3.13a)$$

$$\Lambda(\mu_i) = \Lambda_0(\mu_i) \cup \Lambda_1(\mu_i) := \{\lambda \in CE : Re\,\lambda \geq \mu_i\} \qquad (3.13b)$$

The following result holds concerning an asymptotic comparison of eigensolutions of (21) with the corresponding associated solution of (2.2)-(2.3) under a special form of the perturbation function:

**Theorem 3.2**. Suppose that at least one $\lambda_\mu \in CE$ has real part $\mu$ being identical strict Lyapunov exponent of some solution of (2.2) under (2.3). Suppose also that $f(t, x_t)$ satisfies the hypotheses of Theorem 2.1. Then, a solution of (2.2)-(2.3) satisfies $x(t) = y(t) + O(e^{ct})$, $\forall t \in \mathbf{R}_{0+}$ and any $\mathbf{R} \ni c > \mu$. Also, $x(t) = y(t) + O(t^{\nu-1} e^{\mu t})$ for any $f(t, x_t)$ satisfying the hypotheses of Theorem 2.1 with $\nu$ being the largest multiplicity among those of all distinct $\lambda_\mu \in CE$. Assume, in addition, that $K_{0\alpha} = 0$ for any norm $\alpha$ and $\gamma : \mathbf{R}_{0+} \to \mathbf{R}_{0+}$ satisfies $\gamma(t) = o(e^{-at})$ as $t \to \infty$ for some $a \in \mathbf{R}_+$ (i.e. $\gamma(t) \to 0$ as $t \to \infty$ exponentially). Then, $\exists \varepsilon \in \mathbf{R}_+$ and a nontrivial eigensolution of (2.1) corresponding to the set $\Lambda_0(\mu)$ such that $x(t) = y(t) + o(e^{(\mu-\varepsilon)t})$ as $t \to \infty$.

*Proof:* Note that for any norm $\alpha$, $x_t = O\left(\left|x_t^{P_0(\mathbf{R}_{-h})}\right|_\alpha\right)$; $\forall t \in \mathbf{R}_{0+}$ [Theorem 2.3 (ii)] with $C_e(\mathbf{R}_{-h}) = P_\Lambda(\mathbf{R}_{-h}) \oplus Q_{\Lambda e}(\mathbf{R}_{-h})$, $P_\Lambda(\mathbf{R}_{-h}) = P_{\Lambda_0}(\mathbf{R}_{-h}) \oplus P_{\Lambda_1}(\mathbf{R}_{-h})$ with $P_{\Lambda_0}(\mathbf{R}_{-h})$ being the eigenspace associated with $\Lambda_0(\mu)$ for $t \in \mathbf{R}_{-h}$. For $t \geq h$, the related direct sum decomposition $C(\mathbf{R}_{0+}) = P_\Lambda(\mathbf{R}_{0+}) \oplus Q_\Lambda(\mathbf{R}_{0+})$, $P_\Lambda(\mathbf{R}_{0+}) = P_{\Lambda_0}(\mathbf{R}_{0+}) \oplus P_{\Lambda_1}(\mathbf{R}_{0+})$ since may be used for the eigensolutions of (2.1) since the solution is time-differentiable for $t \in \mathbf{R}_{0+}$. Thus, $x(t) = O(e^{ct})$ and $y(t) = O(e^{ct})$ so that $x(t) - y(t) = O(e^{ct})$ for any given perturbation function $f(t, x_t)$ satisfying the constraints of Theorem 2.1 and all $\mathbf{R} \ni c > \mu = \operatorname{Re} \lambda_\mu$, $\lambda_\mu \in CE$.

If all $\lambda_\mu \in CE$ are simple then $x(t) - y(t) = O(e^{ct})$, $\mathbf{R} \ni c \geq \mu = \operatorname{Re} \lambda_\mu$ for any given perturbation function $f(t, x_t)$ satisfying the constraints of Theorem 2.1.

Otherwise (i.e. at least one $\lambda_\mu \in CE$ is not simple), $x(t) - y(t) = O\left(p_{\lambda_u}(t) e^{\mu t}\right) = O(e^{ct})$ for some $\mathbf{R}_- \ni c > \mu$ with $p_{\lambda_u} : [0, t] \to C^n$ being a polynomial of degree $\nu$ equating the largest multiplicity among those of all distinct $\lambda_\mu \in CE$.

The first part of the result has been fully proved. Now, note that Eq. (3.8) in Lemma 3.1 can be equivalently rewritten as

$$x_t - y_t = \int_0^t \left( T_s(t,\tau) X_0^{P_A} + \left( d[K(t,\tau)]^{Q_A} \right) \right)(\theta) f(\tau, x_\tau) d\tau \,, \quad \mathbf{R}_+ \ni t \geq h \qquad (3.14)$$

Also, if $\gamma(t) = o(e^{-at})$ as $t \to \infty$ and since $|x_t| = o\left(e^{(\mu + \varepsilon) t}\right)$ as $t \to \infty$, irrespective of the multiplicity of $\lambda_\mu$, from the definition of the strict Lyapunov exponent, one gets by using $\gamma(\tau) |x_\tau| \leq M_\gamma M_x e^{(\mu + \varepsilon - \alpha)\tau}$ from (2.7) with $K_{0\alpha} = 0$ (since $\gamma(t) = o(e^{-at}) \to 0$ as $t \to \infty$) and (3.10b):

$$\left| \int_0^t \left( d[K(t,\tau)]^{Q_A} \right) (\theta) f(\tau, x_\tau) d\tau \right|_\alpha \leq M_2(\alpha) e^{(\mu - \varepsilon) t} \int_0^\infty e^{(2\varepsilon - \alpha)\tau} d\tau \leq \frac{M_2(\alpha)}{\alpha - 2\varepsilon} e^{(\mu - \varepsilon) t} = o\left(e^{(\mu - \varepsilon) t}\right)$$
$$, \quad \forall t \in \mathbf{R}_{0+} \qquad (3.15)$$

$$\left| \int_0^t T_s(t,\tau) X_0^{P_A}(\theta) f(\tau, x_\tau) d\tau \right|_\alpha \leq M_3(\alpha) e^{(\mu - \varepsilon) t} \int_0^\infty e^{(2\varepsilon - \alpha)\tau} d\tau \leq \frac{M_2(\alpha)}{\alpha - 2\varepsilon} e^{(\mu - \varepsilon) t} = o\left(e^{(\mu - \varepsilon) t}\right)$$
$$, \quad \forall t \in \mathbf{R}_{0+} \qquad (3.16)$$

since $\left\| X_0^{P_A}(\theta) \right\|_\alpha = \left\| \Phi_\Lambda \Psi_\Lambda(\theta) \right\|_\alpha = \left\| \Phi_\Lambda(0) \Psi_\Lambda(0) e^{2B_\Lambda \theta} \right\|_\alpha < \infty$, $\forall \theta \in [-h, 0]$ provided that $0 < \varepsilon < \frac{\alpha}{2}$. Then, $x_t - y_t = o\left(e^{(\mu - \varepsilon)t}\right)$ as $t \to \infty$, $\forall \varepsilon \in \mathbf{R}_+$ satisfying $0 < \varepsilon < \frac{\alpha}{2}$. □

**Remark 3.3.** An important observation follows. Note that the error between any eigensolutions of the limiting equation (2.1) and the associated solution of (2.2)-(2.3) satisfying $x_t - y_t = o\left(e^{(\mu - \varepsilon)t}\right)$ as $t \to \infty$ if $\gamma(t) = o(e^{-at})$ as $t \to \infty$ and $K_{0\alpha} = 0$ for sufficiently small $\varepsilon \in \mathbf{R}_+$ (proven in Theorem 3.2) does not imply that $x_t = o\left(e^{(\mu - \varepsilon)t}\right)$, $y_t = o\left(e^{(\mu - \varepsilon)t}\right)$ as $t \to \infty$ if $\gamma(t) = o(e^{-at})$ as $t \to \infty$ for sufficiently small $\varepsilon \in \mathbf{R}_+$.

The asymptotic behaviors if $\gamma(t) = o(e^{-at})$ and $K_{0\alpha} = 0$ as $t \to \infty$ are as follows:

$$y_t = o\left(e^{\mu t}\right), \quad x_t = y_t + o\left(e^{(\mu - \varepsilon)t}\right) = o\left(e^{\mu t}\right) + o\left(e^{(\mu - \varepsilon)t}\right) = o\left(e^{\mu t}\right) \text{ as } t \to \infty$$

for sufficiently small $\varepsilon \in R_+$ if $\lambda_\mu \in CE$ is real simple or there are two $\lambda_{\mu_{1,2}} \in CE$ simple complex conjugate ones. In the above , equations, $y_t = o(e^{\mu t})$, $x_t = o(e^{\mu t})$ hold even if $K_{0\alpha} \neq 0$ provided that $\lambda_\mu$ fulfils some of the above constraints. Furthermore,

$$y_t = o(e^{ct}), \quad x_t = y_t + o(e^{(\mu-\varepsilon)t}) = o(e^{ct}) + o(e^{(\mu-\varepsilon)t}) = o(e^{ct}) \text{ as } t \to \infty$$

for sufficiently small $\varepsilon \in R_+$ and all real $c > \mu = Re\,\lambda_\mu$ if $\lambda_\mu \in CE$ is multiple (i.e. of multiplicity greater then unity). Also, $x_t - y_t = O(e^{ct}) \Rightarrow x_t = O(e^{ct}) \wedge y_t = O(e^{ct})$, $\forall t \in R_{0+}$, $\forall c > \mu = Re\,\lambda_\mu$ irrespective of the multiplicity of $\lambda_\mu$ for any generic perturbation function $f(t, x_t)$ satisfying the hypothesis of Theorem 2.1. However, the above implication is true for any arbitrary such a function and $c = \mu$, only if $\lambda_\mu$ is an eigenvalue of multiplicity unity of the limiting equation (2.1). □

Theorem 3.2 leads directly to the subsequent stability result by taking also into account Remark 3.3:

**Corollary 3.4**. Suppose that there exists a (nonnecessarily unique) $\lambda_\mu \in CE$ with $\mu < 0$ being the strict Lyapunov exponent of some solution of (2.2) under (2.3) and that $\Lambda_0(\mu) = \Lambda(\mu)$; i.e. $\Lambda_1(\mu) = \varnothing$. Suppose also that $f(t, x_t)$ satisfies the hypotheses of Theorem 2.1 and, furthermore, $K_{0\alpha} = 0$ (for any $\alpha$-norm) and $\gamma(t) = o(e^{-at})$ as $t \to \infty$ for some $a \in R_+$. Then, $\exists \varepsilon \in R_+$ such that any solutions of (2.1) and (2.2)-(2.3) satisfy:

$$x(t) = y(t) + o\left(e^{-(|\mu|+\varepsilon)t}\right) \to y(t) \to o\left(t^{\nu-1} e^{-|\mu|t}\right) \to o\left(e^{-|c|t}\right) \to 0 \text{ exponentially as } t \to \infty.$$

for some $R_- \ni c \geq \mu$ with $\nu$ being equal to the largest multiplicity of among those of all distinct $\lambda_\mu \in CE$. The above relations hold with $c = \mu$ if and only if all the eigenvalues of (2.1) satisfying the given assumptions are simple. As a result, both functional equations (2.2)-(2.3) and its limiting one are globally asymptotically Lyapunov´s stable with exponential stability.

*Proof*: Theorem 3.2 applies for $c > \mu$ with $|\lambda_\mu|$ being the spectral radius and $\mu < 0$ the spectral (or stability) abscissa, i.e. there is no member of CE with real part to the right of $\mu < 0$ since $\Lambda_1(\mu) = \varnothing$. If all such $\lambda_\mu \in CE$ are simple then the result applies also for $c \geq \mu$ (see Remark 3.2). Define the bounded real non-negative function $V : [t_0, \infty) \times C^n \to R_{0+}$ as $V(t, x(t)) := \|x(t)\|_2^2$ for any finite $t_0 \in R_{0+}$. Note that $V(t, 0) = 0$ and $V(t, x(t))$ is strictly monotonically increasing with $\|x(t)\|$. It follows from (2.2)-(2.3) that $\exists \dot{V} \in PC^{(1)}\left([t_0, \infty) \times C^n, R\right)$, pointwise defined by $\dot{V}(t, x(t)) := 2\dot{x}^*(t) x^*(t)$, and

$0 \leq V(t, x(t)) = V(t_0, x(t_0)) + \int_{t_0}^{t} \dot{V}(\tau, x(\tau)) d\tau \leq K e^{-2|c|(t - t_0)} \to 0$ exponentially as $t \to \infty$.

Then, $\exists \lim_{t \to \infty} \int_{t_0}^{t} \dot{V}(\tau, x(\tau)) d\tau = 2 \lim_{t \to \infty} \int_{t_0}^{t} \dot{x}^*(\tau) x(\tau) d\tau = -V(t_0, x(t_0)) \leq 0$ which is finite since $t_0$ is finite and the solution to (2.2)-(2.3) is everywhere continuous on its definition domain. By assuming $\neg \exists t_0 \in \mathbf{R}_{0+}$ such that $\dot{V}(t, x(t)) \leq 0$, $\forall t \in [t_0, \infty)$ then $\neg \exists \lim_{t \to \infty} \int_{t_0}^{t} \dot{V}(\tau, x(\tau)) d\tau = -V(t_0, x(t_0)) \leq 0$ what is a contradiction to the already proven existence of such a limit. As a result, $\exists t_0$ (sufficiently large but finite) $\in \mathbf{R}_{0+}$ such that $\dot{V}(t, x(t)) \leq 0$, $\forall t \in [t_0, \infty)$ and $\dot{V}(t, x(t)) < 0$ on some subinterval (non-necessarily connected) on infinite measure of $[t_0, \infty)$ with $\dot{V}(t, x(t)) \to 0$ as $t \to \infty$ then necessarily there is a connected terminal subinterval $[t_0', \infty) \subset [t_0, \infty)$ of infinite measure where $\dot{V}(t, x(t)) < 0$. Thus, $V(t, x(t))$ is nonnegative and converges exponentially to zero with non-positive time-derivative which also converges to zero within some interval of infinite measure. As a result, $V(t, x(t))$ is a Lyapunov function with negative time-derivative within a connected real interval of infinite measure which has zero limit. □

Note that the limiting differential functional equation is, furthermore, globally uniformly asymptotically Lyapunov's stable under the asymptotic stability conditions of Corollary 3.4. Note that $K_{0\alpha} > 0$ (for any $\alpha$-norm) in (2.3) implies that $K_{0\alpha'} > 0$ for any other norm $\alpha'$ and conversely. Also, $K_{0\alpha} = 0 \Leftrightarrow K_{0\alpha'} = 0$. The above result concerns with global asymptotic stability with exponential decay rate. Global uniform Lyapunov's stability (i.e. boundedness of solutions with a common upper-bound for all time for any bounded function of initial conditions) holds under weaker conditions; i.e. $K_{0\alpha} > 0$ and $\mu = 0$ if $\lambda_\mu \in CE$ associated with the strict Lyapunov exponent have unity multiplicities or if they have any multiplicities but $K_{0\alpha} = 0$. The precise related stability result follows which proofs follows directly from Corollary 3.4 and Theorem 3.2.

**Corollary 3.5**. Suppose that there exists at least one $\lambda_\mu \in CE$ with $\mu \leq 0$ being the strict Lyapunov exponent of some solution of (2.2) under (2.3) and that $\Lambda_0(\mu) = \Lambda(\mu)$; i.e. $\Lambda_1(\mu) = \varnothing$. Suppose also that $f(t, x_t)$ satisfies the hypotheses of Theorem 2.1 with $K_{0\alpha} \geq 0$ (for any $\alpha$-norm) and $\gamma(t) = o(e^{-at})$ as $t \to \infty$ for some $a \in \mathbf{R}_+$. Then, $\exists \varepsilon \in \mathbf{R}_+$ such that any solutions of (2.1) and (2.2)-(2.3) satisfy:

$$x(t) = y(t) + O\left(e^{-|\mu|t}\right) \to y(t) \to O\left(e^{-|\mu|t}\right), \quad \forall t \in \mathbf{R}_{0+}$$

As a result, the limiting equation (2.1) as well as (2.2)-(2.3) are both globally uniformly Lyapunov's stable if all $\lambda_\mu \in CE$ have any multiplicities and $\mu < 0$ or if they have all unity multiplicities and $\mu = 0$. Eqs (2.1) and (2.2)-(2.3) are both globally Lyapunov´s asymptotically stable with exponential decay rate if $\mu < 0$ and $K_{0\alpha} = 0$, satisfying:

$$x(t) = y(t) + o\left(e^{-|\mu|t}\right) \to y(t) \to o\left(p_{\lambda_u}(t)e^{-|\mu|t}\right) \to o\left(e^{-|c|t}\right) \to 0 \text{ exponentially as } t \to \infty.$$

for some $R_- \ni c \geq \mu$ with $p_{\lambda_u}:[0,t] \to C^n$ being a polynomial of degree equal to the largest multiplicity of all such distinct $\lambda_\mu \in CE$. The above relations hold with $c = \mu$ if and only if all the eigenvalues of (2.1) satisfying the given assumptions are simple, [16]. □

A parallel result to Theorem 3.2 obtained under weaker conditions on the perturbation function (2.3) follows.

**Theorem 3.6**. Suppose that at least one $\lambda_\mu \in CE$ has real part $\mu$ being identical strict Lyapunov exponent of some solution of (2.2) under (2.3). Suppose also that $f(t, x_t)$ satisfies the hypotheses of Theorem 2.1. Suppose, in addition, that $K_{0\alpha} = 0$ for any norm $\alpha$, and $\gamma: R_{0+} \to R_{0+}$ satisfies either $\int_0^\infty t^{\nu-1} \gamma(t) dt < \infty$ or $\gamma(t) \leq o\left(t^{1-\nu}\right)$ as $t \to \infty$, with the inequality being strict if $\nu = 1$, where $\nu$ is the largest multiplicity among those of all distinct $\lambda_\mu \in CE$. Then, $\exists \varepsilon \in R_+$ and a nontrivial eigensolution of (2.1) corresponding to the set $\Lambda_0(\mu)$ such that $x(t) = y(t) + o\left(e^{\mu t}\right)$ as $t \to \infty$. Furthermore, $x(t) = o\left(t^{\nu-1}e^{\mu t}\right)$ and $y(t) = o\left(t^{\nu-1}e^{\mu t}\right)$.

*Proof:* One gets for sufficiently small $\varepsilon \in R_+$ from $\int_0^t t^{\nu-1} \gamma(t) dt < \infty$ and Theorem 3.2, Eqs. (3.15)-(3.16) and $K_{0\alpha} = 0$:

$$\left|\int_0^t \left(d[K(t,\tau)]^{Q_\Lambda}\right)(\theta) f(\tau, x_\tau) d\tau\right|_\alpha \leq M_4(\alpha) e^{\mu t} \int_0^\infty \frac{e^{-\mu\tau} e^{(\varepsilon-\mu)\tau} e^{(\varepsilon+\mu)\tau}}{\tau^{\frac{\nu-1}{2}}} \left(\tau^{\frac{\nu-1}{2}} \gamma(\tau)\right) d\tau$$

(3.17)

$$\leq M_4(\alpha) e^{\mu t} \left[\int_0^\infty \left(\frac{e^{2(2\varepsilon-\mu)\tau}}{\tau^{\nu-1}}\right) d\tau\right]^{1/2} \left[\int_0^\infty \tau^{\nu-1} \gamma^2(\tau) d\tau\right]^{1/2} \leq M_5(\alpha) e^{\mu t} \left[\int_0^\infty \left(\frac{e^{2(2\varepsilon-\mu)\tau}}{\tau^{\nu-1}}\right) d\tau\right]^{1/2}$$

$, \forall t \in R_{0+}$  (3.18)

where $M_5(\alpha) \geq M_4(\alpha) \underset{0 \leq t < \infty}{\text{Sup}} (\gamma(t))$ which is a finite real constant since $\gamma(t)$ is continuous on $R_{0+}$ and has zero limit as $t \to \infty$ so that it is uniformly bounded. Now, if $\mu \in R_+$ for $R_+ \ni \varepsilon \in \left(0, \frac{\mu}{2}\right), R_- \ni d > b := (2\varepsilon - \mu)\nu \in R_-$ and $R_{0+} \ni t \geq t_0 := e^{\frac{d-b}{\nu-1}}$:

$$\int_{t_0}^\infty \left(\frac{e^{2(2\varepsilon-\mu)\tau}}{\tau^{\nu-1}}\right) d\tau \leq \int_{t_0}^\infty e^{-|d|\tau} d\tau \leq \frac{1}{|d|} = (\mu - 2\varepsilon)\nu < \infty; \int_0^{t_0} \left(\frac{e^{2(2\varepsilon-\mu)\tau}}{\tau^{\nu-1}}\right) d\tau \leq M_6 < \infty$$

(3.19)

since $t_0$ is finite. Now from (3.19) into (3.18):

$$\left|\int_0^t \left(d[K(t,\tau)]^{Q_A}\right)(\theta)f(\tau,x_\tau)d\tau\right|_\alpha \leq M_5(\alpha)\left(\sqrt{M_6} + \sqrt{(\mu-2\varepsilon)v}\right)e^{\mu t} = o\left(e^{\mu t}\right)$$

as $t \to \infty$ \quad (3.20)

if $\mu \in \mathbf{R}_+$. If $\mu \in \mathbf{R}_{0-}$ then Eqs. (3.10) of Lemma 3.1 may be replaced for $\mu \in \mathbf{R}_{0-}$ by taking $\varepsilon = 0$ with

$$\left|T(t,0)\phi^{P_A}\right| \leq M_1 t^{\nu-1} e^{-|\mu|t}\left|\phi^{P_A}\right|, \quad \left|T(t,0)X_0^{P_A}\right| \leq M_1 t^{\nu-1} e^{-|\mu|t}; \quad \forall t \in \mathbf{R}_{0-} \quad (3.21)$$

$$\left|T(t,0)\phi^{Q_{Ae}}\right| \leq M_1 t^{\nu-1} e^{-|\mu|t}\left|\phi^{Q_{Ae}}\right|; \quad \forall t \in \mathbf{R}_{0+} \quad (3.22)$$

By using (3.22), Eq. (3.18) may be replaced $\mu \in \mathbf{R}_{0-}$ with

$$\left|\int_0^t \left(d[K(t,\tau)]^{Q_A}\right)(\theta)f(\tau,x_\tau)d\tau\right|_\alpha \leq M_4(\alpha)\int_0^t e^{-|\mu|(t-\tau)}\left(\tau^{\nu-1}\gamma(\tau)\right)d\tau$$

$$\leq M_4(\alpha)\int_0^\infty \tau^{\nu-1}\gamma(\tau)d\tau \leq M_7(\alpha) < \infty, \quad \forall t \in \mathbf{R}_{0+}$$

On the other hand, $\int_0^\infty t^{\nu-1}\gamma(\tau)d\tau < \infty \Rightarrow \gamma(t) = o\left(\dfrac{1}{t^{\nu-1}}\right)$ as $t \to \infty$ so that one gets for $\mu \in \mathbf{R}_{0-}$ from Theorem 3.2, Eq. (3.15), since the integrals in (3.19) are bounded:

$$\left|\int_0^t \left(d[K(t,\tau)]^{Q_A}\right)(\theta)f(\tau,x_\tau)d\tau\right|_\alpha \leq M_4(\alpha)e^{\mu t}\int_0^t \tau^{1-\nu} e^{\varepsilon(2\tau-t)}d\tau$$

$$\leq 2M_4(\alpha)e^{\mu t}\int_0^{2t}\left(\dfrac{\tau}{2}\right)^{1-\nu} e^{\varepsilon(\tau-t)}d\tau \leq 2M_4(\alpha)e^{\mu t}\int_0^\infty \left(\dfrac{\tau}{2}\right)^{1-\nu} e^{\varepsilon(\tau-t)}d\tau$$

$$\leq M_8(\alpha)e^{\mu t} = o\left(e^{-|\mu|t}\right) \text{ as } t \to \infty, \text{ irrespective of the norm } \alpha.$$

As a final result, one has $\left|\int_0^t \left(d[K(t,\tau)]^{Q_A}\right)(\theta)f(\tau,x_\tau)d\tau\right|_\alpha = o\left(e^{\mu t}\right)$ as $t \to \infty$, $\forall \mu \in \mathbf{R}$, and, furthermore, $\left|\int_0^t \left(d[K(t,\tau)]^{Q_A}\right)(\theta)f(\tau,x_\tau)d\tau\right|_\alpha = O\left(e^{\mu t}\right), \forall \mu \in \mathbf{R}_{0-}$, $\forall t \in \mathbf{R}_{0+}$, irrespective of the norm $\alpha$. Similar properties follow under close proofs for $\left|\int_0^t T_s(t,\tau)X_0^{P_A}(\theta)f(\tau,x_\tau)d\tau\right|_\alpha$. Thus, $x_t - y_t = o\left(e^{\mu t}\right)$ as $t \to \infty$ by using the above results in (3.18) and $x_t - y_t = O\left(e^{\mu t}\right), \forall t \in \mathbf{R}_{0+}$. The properties $x(t) = o\left(t^{\nu-1}e^{\mu t}\right)$ and $y(t) = o\left(t^{\nu-1}e^{\mu t}\right)$ of the solutions of (2.1) and (2.2)-(2.3) under the given hypothesis follow directly by using the appropriate reasoning quoted from Remark 3.2. □

Lyapunov's stability properties in terms of boundedness of the solutions and their asymptotic or exponential convergence to the equilibrium obtained from Theorem 3.6 are immediate as given in the following direct result.

**Corollary 3.7**. Corollaries 3.4-3.5 also apply "mutatis-mutandis" under the assumptions of Theorem 3.6.

□

## 4. Some direct consequences and applications

Some particular cases of the functional differential equations (2.1) and (2.2)-(2.3) are of interest concerning stability issues as follows:

*4.1 Functional differential equations with point delays*

The functional differential equation (2.2) can be equivalently described in the absence of finite distributed delays and Volterra- type dynamics by the n-the order system of n functional first-order differential equations:

$$\dot{x}(t) = \sum_{i=0}^{m} A_i x(t - h_i) + f(t, x_t) \tag{4.1}$$

where $h_0 = 0$. Its limiting differential system is defined in the same way from (4.1) for a perturbation function $f(t, x_t) \equiv 0$. The generic form of such a function is given in (2.2) and can include nonzero dynamics of finite distributed delays and Volterra- type dynamics. Thus, $\hat{T}(s) = \left( sI_n - \sum_{i=0}^{m} A_i e^{-h_i s} \right)^{-1}$ everywhere the inverse exists. Remember for later use that a matrix-valued function $G: C_{0+} \to C^{r \times m}$ is in the Hardy space $H_\infty$ if it is analytic in $C_+$, $\lim_{\sigma \to 0^+} G(\sigma + i\omega) = G(i\omega)$ for almost all $\omega \in R_{0+}$ and $\|G\|_\infty := \sup_{s \in C_{0+}} \bar{\sigma}(G(s)) = \sup_{\omega \in R} \bar{\sigma}(G(i\omega)) < \infty$ with $\bar{\sigma}(G)$ denoting the largest singular value of $G$. $RH_\infty$ is a the subset of $H_\infty$ of real-rational matrix valued functions then being proper and stable (i.e. they have no more zeros than poles and all the poles are in $C_-$), [21]. Note that the used norm notation $\|G\|_\infty$ for $G \in H_\infty$ is similar to the notation for $\ell_\infty$- matrix /vector norms but no confusion is expected from the different context of use. If $G(s) = G^*(s^*)$ then $\|G\|_\infty := \sup_{\omega \in R_{0+}} \bar{\sigma}(G(i\omega))$

The following result holds:

**Theorem 4.1**. The following properties hold:
(i) The limiting differential system associated with (4.1) is globally asymptotically stable independent of the delays (i.e. $\forall h_i \in R_{0+}$ ($\forall i \in \overline{m}$) if $A_0$ is a stability matrix and there exist $\beta_i \in R_+$ ($\forall i \in \overline{m}$) such that $\sum_{i=1}^{m} \beta_i^2 = 1$, and

$$\left\| \left( i\omega I_n - \sum_{i=0}^{m} A_i \right)^{-1} \left[ \beta_1^{-1} A_1, \ldots, \beta_r^{-1} A_m \right] \right\|_2 < 1, \forall \omega \in R_{0+} \tag{4.2}$$

or

$$\left\| \left[ \beta_1^{-1} A_1, \ldots, \beta_r^{-1} A_m \right] \right\|_2 < -\kappa_2(A_0) \tag{4.3}$$

This implies as a result that $\kappa_2(A_0) < 0$ and that any matrix $\left(A_0 + \sum_{i \in J_k} A_i\right)$ is a stability matrix, $\forall J_k \subseteq \overline{m}$, where $J_{\overline{m}} := \bigcup_{k \in s_{\overline{m}}} J_k$ is the union of all the $s_{\overline{m}}$ denumerable sets obtained by combining members of $\overline{m}$ (i.e. $J_{\overline{m}}$ is the set of parts of $\overline{m}$).

(ii) If Property (i) holds then $\overline{A}_0 = A_0 + \rho I_n$ is a stability matrix for all $R_{0+} \ni \rho \in [0, \rho_0)$ with $\rho_0 := 1/\|A_0^{-1}\|_2$. The differential system (4.1) is globally asymptotically stable independent of the size of the delays with stability (or spectral) abscissa of at least $-\rho_0 := -\min(\rho_{01}, \rho_{02}) < 0$, with

$$\rho_{01} := 1/\|A_0^{-1}\|_2 \tag{4.4}$$

$$\rho_{02} := 1/\left\|\left(sI_n - \sum_{i=0}^m A_i\right)^{-1}\right\|_\infty = \sup_{s \in C_{0+}} \overline{\sigma}\left(\left(sI_n - \sum_{i=0}^m A_i\right)^{-1}\right) = \sup_{\omega \in R_{0+}} \overline{\sigma}\left(\left(i\omega I_n - \sum_{i=0}^m A_i\right)^{-1}\right)$$

(4.5)

provided that:

. $A_0$ is a stability matrix

. Eq. (4.2) holds for some $\beta_i \in R_+$ ($\forall i \in \overline{m}$) such that $\sum_{i=1}^m \beta_i^2 = 1$

or if $\rho \in [0, \kappa_2(A_0))$ and, furthermore,

$$\left\|\left[\beta_1^{-1} A_1, \ldots, \beta_r^{-1} A_m\right]\right\|_2 < -\kappa_2(A_0) + \rho < 0 \tag{4.6}$$

so that the limiting differential system (4.1) is globally asymptotically stable independent of the delays with stability abscissa of at least $-\kappa_2(A_0) < 0$

*Proof*: (i) Property (i) follows from results in [20, Proposition 4.8] with Eq. (4.2) as follows. Note that $h_i = 0$, $\forall i \in \overline{m}$, the (delay-free) system (4.1) is globally asymptotically stable so that $\left\|\left(sI_n - \sum_{i=0}^m A_i\right)^{-1}\right\|_\infty < \infty$ so that $\left\|\left(sI_n - \sum_{i=0}^m A_i e^{-h_i s}\right)^{-1}\right\|_\infty < \infty$. Then, the limiting differential system associated with (4.1) is globally asymptotically stable independent of the delays if Eq. (4.2) holds which is obtained by guaranteeing the inverse below on $C_{0+}$ under a necessary and sufficiency- type text on the complex imaginary axis:

$$\left(sI_n - \sum_{i=0}^m A_i e^{-h_i s}\right)^{-1} = \left(\sum_{i=0}^m A_i\left(1 - e^{-h_i s}\right)\right)^{-1}\left(sI_n - \sum_{i=0}^m A_i\right)^{-1}$$

(ii) Since property (i) holds, then Eq. (4.2) holds so that there exist $\beta_i \in \mathbf{R}_+$ $(\forall i \in \overline{m})$ (which are in general distinct from those in Property (i) but we keep the same notation) such that $\sum_{i=1}^{m} \beta_i^2 = 1$, and for any $\mathbf{R}_{0+} \ni \rho \in [0, \rho_{02})$, $\forall \omega \in \mathbf{R}_{0+}$:

$$\left\| \left( sI_n - \sum_{i=0}^{m} A_i \right) \left( I_n - \rho \left( sI_n - \sum_{i=0}^{m} A_i \right)^{-1} \right)^{-1} \left[ \beta_1^{-1} A_1, \ldots, \beta_r^{-1} A_m \right] \right\|_\infty$$

$$= \left\| \left( i\omega I_n - \sum_{i=0}^{m} A_i \right) \left( I_n - \rho \left( i\omega I_n - \sum_{i=0}^{m} A_i \right)^{-1} \right)^{-1} \left[ \beta_1^{-1} A_1, \ldots, \beta_r^{-1} A_m \right] \right\|_2$$

$$\leq \left\| \left( i\omega I_n - \sum_{i=0}^{m} A_i \right)^{-1} \left[ \beta_1^{-1} A_1, \ldots, \beta_r^{-1} A_m \right] \right\|_2 \left( \frac{1}{1 - \rho \left\| \left( i\omega I_n - \sum_{i=0}^{m} A_i \right)^{-1} \right\|_2} \right) < 1$$

(4.7)

and the above $H_\infty$-norm exists if $\left\| \left( sI_n - \sum_{i=0}^{m} A_i \right)^{-1} \right\|_\infty < \infty$ from Banach´s Perturbation Lemma, [22]. On the other hand, if $A_0$ is stability matrix then $A_0 + \rho I_n$ is still a stability matrix if Eq. (4.4) holds by applying the min max computation approach for the eigenvalues of $A_0 + \rho I_n$ which are largest than those of $A_0$ since the identity matrix is positive definite. As a result, if global asymptotic stability independent of delays holds then it also holds with stability abscissa $-\rho_0 < 0$. The first part of property (ii) has been proven. The second part follows by replacing Eq. (4.3), through the use of the properties of the matrix measure [20], [22], with

$$\left\| \left[ \beta_1^{-1} A_1, \ldots, \beta_r^{-1} A_m \right] \right\|_2 < -\kappa_2 (A_0 + \rho I_n) \leq -\kappa_2 (A_0) + \rho < 0 \qquad \square$$

Theorem 4.1 can be directly combined with Corollaries 3.4-3.5 as follows:

**Corollary 4.2**. Assume that the functional limiting differential system (4.1) satisfies Theorem 4.1. Thus, both the perturbed and the nominal functional differential systems have a negative strict Lyapunov exponent which satisfies $\mu \leq -min(\rho_{01}, \rho_{02}, \kappa_2(A_0)) < 0$. The solution of (4.1), subject to (2.3) with the hypotheses in Theorem 2.1, and that of its limiting differential system are both globally asymptotically stable independent of the delays.

If, furthermore, $\gamma: \mathbf{R}_{0+} \to \mathbf{R}_{0+}$ satisfies one of the subsequent conditions:

(a) $\gamma(t) = o(e^{-at})$ as $t \to \infty$

(b) $\int_0^\infty t^{\nu-1} \gamma(t) dt < \infty$ as $t \to \infty$

(c) $\gamma(t) \leq o(t^{1-\nu})$ as $t \to \infty$

then any solution of the perturbed functional differential system (4.1) fulfils either Corollaries 3.4- 3.5 to Theorem 3.2 [under the condition (a)] or Corollary 3.7 to Theorem 3.6 [under the condition (b) or the condition (c)] . □

A parallel result to theorem 4.1 (i) may be obtained from the subsequent remark.

**Remark 4.3**. Note tat Eq. (4.2) holds if

$$\left\| (sI_n - A_0)\left(I_n - (sI_n - A_0)^{-1}\left(\sum_{i=0}^{m} A_i\right)\right)\left(I_n - \rho\left(sI_n - \sum_{i=0}^{m} A_i\right)^{-1}\right)^{-1} \left[\beta_1^{-1} A_1, \ldots, \beta_r^{-1} A_m\right] \right\|_\infty < 1$$

, $\forall s \in C_{0+}$ provided that $\left((i\omega I_n - A_0)\left(I_n - (i\omega I_n - A_0)^{-1}\left(\sum_{i=0}^{m} A_i\right)\right)\right)^{-1}$ exists , $\forall \omega \in R_{0+}$.

This inverse exists provided that $A_0$ is a stability matrix, as it has been assumed, since then it has no imaginary eigenvalues or real ones at the origin and, furthermore, $\left\|(sI_n - A_0)^{-1}\left(\sum_{i=1}^{m} A_i\right)\right\|_\infty \leq a_0 < 1$. If

$$\rho \left\|\left(sI_n - \sum_{i=0}^{m} A_i\right)^{-1}\right\|_\infty \leq \rho \left\|(sI_n - A_0)^{-1}\left(\sum_{i=1}^{m} A_i\right)\right\|_\infty \leq \rho a_0 \leq b_0 < 1$$

$$\left\| \left((i\omega I_n - A_0)\left(I_n - (i\omega I_n - A_0)^{-1}\left(\sum_{i=0}^{m} A_i\right)\right)\right)^{-1} \left[\beta_1^{-1} A_1, \ldots, \beta_r^{-1} A_m\right] \right\|_2 < 1, \quad \forall \omega \in R_{0+}$$

provided that $\left((i\omega I_n - A_0)\left(I_n - (i\omega I_n - A_0)^{-1}\left(\sum_{i=0}^{m} A_i\right)\right)\right)^{-1}$ exists , $\forall \omega \in R_{0+}$. This inverse exists provided that $A_0$ is a stability matrix, as it has been assumed, since then it has no imaginary eigenvalues or real ones at the origin and furthermore $\left\|(sI_n - A_0)^{-1}\left(\sum_{i=1}^{m} A_i\right)\right\|_\infty \leq a_0 < 1$. Then, the system (4.1) is globally asymptotically stable independent of the delays if $\left\|(sI_n - A_0)^{-1}\left(\sum_{i=1}^{m} A_i\right)\right\|_\infty \leq a_0 < 1$ (what implies that $A_0$ is a stability matrix) and $\left\|[\beta_1^{-1} A_1, \ldots, \beta_r^{-1} A_m]\right\|_2 < a_0^{-1}$. The more general constraint Eq. (4.7) is guaranteed with any $R_{0+} \ni \rho < a_0^{-1}$ if $\left\|[\beta_1^{-1} A_1, \ldots, \beta_r^{-1} A_m]\right\|_2 < \frac{1 - \rho a_0}{a_0}$ which guarantees global asymptotic stability independent of the delays of the limiting differential system of (4.1) and $A_0$ is a stability matrix of spectral abscissa of at least $-\rho < 0$. □

*4.2 Functional differential equations with mixed finite point and varying distribute delays*

Assume that (2.2)-(2.3) consist of a single finite constant point delay $h_1 \geq 0$ and a single distributed time-varying one $h:[-1,0] \to R_{0+}$ leading to the functional differential system:

$$\dot{x}(t) = \sum_{i=0}^{1} A_i x(t-h_i) + \int_{-1}^{0} d\alpha(\tau) x(t-h(\tau)) + f(t, x_t) \tag{4.8}$$

with $h_0 = 0$. The perturbation function is defined by (2.3) and satisfies the hypothesis of Theorem 2.1. The limiting equation is defined for identically zero perturbation function. Eq. (4.8) is equivalent to:

$$\dot{x}(t) = (A_0 + A_1 + \alpha(0) - \alpha(-1))x(t) + A_1(x(t-h_1) - x(t)) + \int_{-1}^{0} d\alpha(\tau)(x(t-h(\tau)) - x(t)) + f(t, x_t) \tag{4.9}$$

so taht the characteristic equation of the limiting equation is:

$$det\left(sI_n - A_0 - A_1 e^{-h_1 s} - \int_{-1}^{0} \lambda(\tau) d\alpha(\tau)\right) \tag{4.10}$$

$$= det\left(sI_n - (A_0 + \alpha(0) - \alpha(-1)) - A_1 e^{-h_1 s} - \int_{-1}^{0} \lambda(\tau) d\alpha(\tau)\right) \tag{4.11}$$

$$= det\left(sI_n - (A_0 + A_1 + \alpha(0) - \alpha(-1)) - A_1(e^{-h_1 s} - 1) - \int_{-1}^{0} \lambda(\tau) d\alpha(\tau)\right) = 0 \tag{4.12}$$

for any continuous function $\lambda(t)$ mapping $[-1,0]$ on the unit circle centred at the origin of the complex plane provided that $\alpha$ is non-atomic at zero. The following result is concerned with the global asymptotic stability of the differential system (4.8):

**Theorem 4.4**. The following properties hold:

(i) Assume that $A_0$ is a stability matrices, so that $\exists \rho_1 \in R_+ : \kappa_2(A_0) \leq -\rho_1 < 0$. Thus, if

$$\frac{1}{\|A_1\|_2 + \left\|\int_{-1}^{0} \lambda(\tau) d\alpha(\tau)\right\|_2} \left\|(sI_n - A_0 - \rho_1 I_n)^{-1}\right\|_\infty < 1 \tag{4.13}$$

then the limiting equation associated with (4.8) is globally asymptotically stable independent of the delays.

(ii) Assume that $A_0$ and $(A_0 + \alpha(0) - \alpha(-1))$ are both stability matrices, so that $\exists \rho_2 \in R_+ : \kappa_2(A_0 + \alpha(0) - \alpha(-1)) \leq -\rho_2 < 0$. Thus, if

$$\frac{1}{\|A_1\|_2 + \|\alpha(0)\|_2 + \|\alpha(1)\|_2} \left\|(sI_n - A_0 - \alpha(0) + \alpha(-1) - \rho_2 I_n)^{-1}\right\|_\infty < 1 \tag{4.14}$$

then the limiting equation associated with (4.8) is globally asymptotically stable independent of the delays.

(iii) Assume that $A_0 + A_1 + \alpha(0) - \alpha(-1)$ is stability matrix so that $\exists \rho_3 \in \mathbb{R}_+$: $\kappa_2(A_0 + A_1 + \alpha(0) - \alpha(-1)) \leq -\rho_3 < 0$. Thus, if

$$\frac{1}{2\|A_1\|_2 + \|\alpha(0)\|_2 + \|\alpha(1)\|_2} \left\|(sI_n - A_0 - A_1 - \alpha(0) + \alpha(-1) - \rho_3 I_n)^{-1}\right\|_\infty < 1 \qquad (4.15)$$

then the limiting equation associated with (4.8) is globally asymptotically stable independent of the delays.

*Proof*: Property (i) is direct by using a similar reasoning to that in Theorem 4.1 by extending a result in [20, Section 4.4.5] for a single distributed delayed dynamics, the limiting differential system associated with (4.8) is globally asymptotically stable if and only if for any non- characteristic zero of the limiting equation, $det\left(sI_n - A_0 - A_1 e^{-hs} - \int_{-1}^{0} \lambda(\tau) d\alpha\tau\right) \neq 0$, $\forall s \in \mathbb{C}_{0+}$, or equivalently, if the inverse below exists for $\forall s \in \mathbb{C}_{0+}$ by using the re-arrangement Eq. (4.10):

$$\left(sI_n - A_0 - A_1 e^{-hs} - \int_{-1}^{0} \lambda(\tau) d\alpha(\tau)\right)^{-1} = \left((sI_n - A_0 - \rho_1 I_n)\left(I_n - (sI_n - A_0 - \rho_1 I_n)^{-1}\left(A_1 e^{-hs} + \int_{-1}^{0} \lambda(\tau) d\alpha(\tau)\right)\right)\right)^{-1}$$

(4.16)

Such an inverse exists within $\mathbb{C}_{0+}$ under the conditions in Property (i) so that the evolution operator exists $T(t,0)$ since its associate resolvent exists and it is compact in $\mathbb{C}_{0+}$. Properties (ii)-(iii) are direct alternative sufficient conditions to those in Property (i) by using a similar re-arrangements for an inverse matrix as Eq. (4.16) by using Eq. (4.11) and Eq.(4.12), respectively, instead of Eq. (4.10). □

A parallel result to Corollary 4.2 now follows from Theorem 4.4.

**Corollary 4.5**. Assume that the functional limiting differential system of equations (4.8) satisfies any of the Properties of Theorem 4.4. Thus, both the limiting and the perturbed functional differential equation have a negative strict Lyapunov exponent satisfying $\mu \leq -\rho < 0$ with $\rho$ equalizing the corresponding $\rho_i$ (i = 1, 2, 3). Corollary 4.2 holds "mutatis- mutandis". □

The extension of Theorem 4.4 and Corollary 4.5 to the case of multiple point and distributed delays is direct and then omitted.

ACKNOWLEDGMENTS

The author is very grateful to the Spanish Ministry of Education by its partial support of this work through project DPI2006-00714. He is also grateful to the Basque Government by its support through GIC07143-IT-269-07, SAIOTEK SPED06UN10 and SPE07UN04.